\documentclass[12pt,draft]{article}

\usepackage{amsthm,amssymb,amsmath,amscd,amstext,mathrsfs}
\usepackage{latexsym,mathptmx,xypic,tikz,caption}

\newtheorem{theorem}{Theorem}[section]
\newtheorem{corollary}{Corollary}[section]

\newtheorem{lemma}{Lemma}[section]
\newtheorem*{remark}{{\it Remark}}

\newcommand{\nc}{\newcommand} 

\nc{\C}{{\mathbb C}}
\nc{\R}{{\mathbb R}}
\nc{\HH}{{\mathbb H}}
\nc{\Z}{{\mathbb Z}}
\nc{\N}{{\mathbb N}}
\nc{\dd}{{\rm d}}

\nc{\ii}{{\bf i}}
\nc{\cg}{{\mathscr G}}
\nc{\co}{{\mathscr O}}
\nc{\crr}{{\mathscr R}}

\begin{document}

\title{Global solvability of the vacuum Einstein equation and\\
the strong cosmic censorship in four dimensions}

\author{G\'abor Etesi\\
\small{{\it Department of Geometry, Mathematical Institute, Faculty of 
Science,}}\\
\small{{\it Budapest University of Technology and Economics,}}\\
\small{{\it Egry J. u. 1, H \'ep., H-1111 Budapest, Hungary}}
\footnote{E-mail: {\tt etesi@math.bme.hu}}}

\maketitle
                          
\pagestyle{myheadings}

\markright{Solvability of the vacuum Einstein equation in four dimensions}

\thispagestyle{empty}

\begin{abstract} 
Let $M$ be a connected, simply connected, oriented, closed, 
smooth four-manifold which is spin (or equivalently 
having even intersection form) and put $M^\times:=M\setminus\{{\rm point}\}$. 
In this paper we prove that if $X^\times$ is a smooth four-manifold 
homeomorphic but not necessarily diffeomorphic to $M^\times$ 
(more precisely, it carries a smooth structure {\it \`a la} Gompf) then 
$X^\times$ can be equipped with a complete Ricci-flat Riemannian metric. As a 
byproduct of the construction it follows that this metric is self-dual as well 
consequently $X^\times$ with this metric is in fact a hyper-K\"ahler manifold. 
In particular we find that the largest member of the 
Gompf--Taubes radial family of large exotic $\R^4$'s admits a complete 
Ricci-flat metric (and in fact it is a hyper-K\"ahler manifold). 

These Riemannian solutions are then converted into Ricci-flat Lorentzian ones 
thereby exhibiting lot of new vacuum solutions which are not accessable by the 
initial vaule formulation. A natural physical interpretation of them in the 
context of the strong cosmic censorship conjecture and topology change is 
discussed.
 
\end{abstract}

\centerline{AMS Classification: Primary: 83C05, Secondary: 57N13, 53C28}
\centerline{Keywords: {\it Ricci-flatness; Exotic $\R^4$; Twistors; 
Strong Cosmic Censorship}}

%%%%%%%%%%%%%%%%%%%%%%%%%%%%%%

\section{Introduction and summary}
\label{one}

%%%%%%%%%%%%%%%%%%%%%%%%%%%%%%%%%%%%%

Solving the {\it vacuum Einstein equation} globally, or in other words: 
finding a (pseudo-)Riemannian {\it Ricci-flat} metric along a 
differentiable manifold i.e., a metric $g$ which satisfies the second order 
non-linear partial differential equation 
\[{\rm Ric}_g=0\] 
over a differentiable manifold $M$, is a century-old evergreen problem 
dwelling in the heart of modern differential geometry \cite{bes} and 
theoretical physics \cite{wald}. The problem of solvability naturally 
splits up into {\it local} and {\it global} solvability and also depends on 
the signature of the metric. Let us first consider the {\it Riemannian} case. 
Thanks to its non-linearity, solvability of the Ricci-flatness condition is 
already locally problematic; nevertheless exploiting its elliptic character 
various kinds of local existence results (e.g. \cite{det,gas}) are known 
at least for the related but in some sense complementary equation 
${\rm Ric}_g=\Lambda g$ (whose solutions are called {\it Einstein 
metrics}) with $\Lambda\not=0$. As one expects, in these local 
existence problems the dimension of $M$ plays no special role. However 
dimensionality issues characteristically enter the 
game when one considers global solvability. Finding global 
solutions in four dimensions i.e., when $\dim_\R M=4$ is particularly 
important from a physical point of view and quite interestingly, from the 
mathematical viewpoint, precisely this is the dimension where global 
solvability is the most subtle. As it is well-known, if $\dim_\R M<4$ the 
vacuum Einstein equation reduces to a full flatness condition on the metric 
hence it admits only a ``few'' global solutions; on the contrary, if 
$\dim_\R M>4$ there are no (known) obstructions for global solvability hence 
apparently there are ``too many'' global solutions. A delicate balance is 
achieved if $\dim_\R M=4$: for instance by a classical result \cite{hit0,tho} 
we know that a Riemannian Einstein (hence in particular a Ricci-flat) metric 
on a compact $M$ can exist only if 
its Euler characteristic $\chi (M)$ and signature $\sigma (M)$ obey the 
inequality $\chi (M)\geqq\frac{3}{2}\vert\sigma (M)\vert$. This implies 
for example that the connected sum of at least five copies 
of complex projective spaces cannot be Einstein. However even in four 
dimensions if $M$ is non-compact there are no (known) obstruction against 
the solvability of the vacuum Einstein equation.

Restricting attention to the four dimensional case from now on, the main 
result of the paper---strongly motivated by \cite{che-nem} and considered 
as a substantially improved and technically revised and greatly simplified 
version of our earlier efforts \cite{ete1, ete2}---can be formulated in the 
{\it Riemannian} setting as 

\begin{theorem} Let $M$ be a connected, simply connected, oriented, closed 
(i.e., compact without boundary), smooth $4$-manifold which is spin (or 
equivalently having even intersection form) and take the punctured space 
$M^\times:=M\setminus\{{\rm point}\}$. If $X^\times$ is a smooth 
$4$-manifold {\it homeomorphic} but not necessarily diffeomorphic to 
$M^\times$ such that it carries a smooth structure {\it \`a la} Gompf then 
$X^\times$ can be equipped with a complete Ricci-flat Riemannian metric. 
\label{fotetel1} 
\end{theorem}

\noindent As an extreme but important application of Theorem 
\ref{fotetel1} we obtain

\begin{corollary}
Let $R^4$ be the largest member of the Gomp--Taubes radial family of large
exotic $\R^4$'s. Then $R^4$ carries a complete Ricci-flat Riemannian metric.
\label{kovetkezmeny1}
\end{corollary}

\noindent The proof of Theorem \ref{fotetel1} is based on a successive 
application of basic results by Gompf \cite{gom1,gom2,gom3}, Penrose 
\cite{pen1}, Taubes \cite{tau1,tau2} and Uhlenbeck \cite{uhl} on exotic 
smooth structures, twistor theory, self-dual spaces and singularity 
removal in Yang--Mills fields, respectively. The idea in the spirit of 
twistor theory is to convert the real-analytic problem of solving ${\rm 
Ric}_g=0$ on the real $4$-space $M^\times$ into a complex-analytic problem 
on a complex $3$-space $Z$ associated to $M^\times$. This is in principle 
simple and works as follows. Take an arbitrary oriented and closed smooth 
$4$-manifold $M$. In the first step, following Taubes, by connected 
summing sufficiently (but finitely) many complex projective spaces to $M$, 
we construct a space $\overline{X}_M\cong M\#\C P^2\#\dots\#\C P^2$ which 
(with respect to its induced orientation) carries a self-dual metric 
$\overline{\gamma}$. Then, in the second step following Penrose, we 
observe that $\overline{Z}$, the {\it twistor space} of 
$(\overline{X}_M,\overline{\gamma})$, is a complex $3$-manifold. Let 
$X_M\subset\overline{X}_M$ be the open space obtained by deleting 
carefully choosen closed subsets, homeomorphic to a projective line, from 
every $\C P^2$ factor of $\overline{X}_M$ and put 
$\gamma:=\overline{\gamma}\vert_{X_M}$ and $Z:=\overline{Z}\vert_{X_M}$. 
Making use of $Z$ we can conformally rescale the incomplete self-dual 
space $(X_M,\gamma )$ to a complete Ricci-flat one $(X_M,g)$ if $M$ is 
connected, simply connected and spin. In the third and last step, by the 
aid of Uhlenbeck's singularity removal theorem, we remove (or fill in) the 
extra $\C P^2\setminus\C P^1=\R^4$'s along $X_M$ to obtain an open smooth 
space $X^\times$ which is homeomorphic to the punctured space $M^\times$ 
however is not necessarily diffeomorphic to it by results of Gompf. The 
result is a connected, simply connected, open, complete, Ricci-flat 
Riemannian spin $4$-manifold $(X^\times, g)$.

By the conformal invariance of self-duality this technical condition in 
fact survives the whole procedure. Taking into account that a complete 
Ricci-flat and self-dual metric on a simply connected $4$-manifold always 
induces a hyper--K\"ahler structure on it \cite[Chapter 13]{bes}, we can 
re-formulate the result of our construction as

\begin{theorem}
The complete Ricci-flat metric of Theorem \ref{fotetel1} on
$X^\times$ with its fixed orientation is self-dual as well consequently
$X^\times$ carries a hyper-K\"ahler structure, too.
\label{fotetel2}
\end{theorem}

\noindent In this way we obtain

\begin{corollary}
The space $R^4$ of Corollary \ref{kovetkezmeny1} carries a hyper-K\"ahler   
structure.
\label{kovetkezmeny2}
\end{corollary}

Next let consider the analogous problem in {\it Lorentzian} signature. 
Surely the most productive---and both mathematically and physically 
extraordinary important---presently known method to find global solutions of 
the Lorentzian vacuum Einstein equation is based on the {\it initial 
value formulation} \cite[Chapter 10]{wald} which exploits the hyperbolic 
character of the Ricci-flatness condition (far from being complete, just for 
recent results cf. eg. \cite{car-sch,cor,ise-maz-pol, maz-pol-uhl}) and the 
references therein). In this approach one starts with an appropriate 
initial value data set, subject to the (simpler) vacuum constraint equations, 
on a {\it three} dimensional manifold $\Sigma$ and obtains solutions of the 
original vacuum Einstein equation on a {\it four} dimensional manifold $M$ 
which is {\it always diffeomorphic} to the smooth product $\Sigma\times\R$ 
(with the unique smooth structures on the factors) \cite{ber-san, che-nem}. 
It is worth calling attention that even if the initial value formulation 
produces an abundance of solutions from the viewpoint of 
{\it global analysis} and {\it theoretical physics}, it is quite inproductive 
from the viewpoint of (low dimensional) {\it differential topology}. 
To illustrate this, suppose we want to find spaces $(M,g)$ satisfying 
${\rm Ric}_g=0$ over a connected and simply connected, open four-manifold $M$. 
If the initial value formulation is applied, and if in this case 
we impose a further condition that the corresponding Cauchy surface $\Sigma$ be 
compact, then by the Poincar\'e--Hamilton--Perelman theorem $\Sigma$ must 
be homeomorphic hence diffeomorphic to the three-sphere $S^3$ consequently $M$ 
is uniquely fixed to be $S^3\times\R$ up to diffeomorphisms (but of course 
this unique $M$ still can carry lot of non-isometric Ricci-flat metrics $g$). 

However in sharp contrast to this differentio-topological 
rigidity of initial value formulation in the simply connected setting we 
obtain 

\begin{theorem} Consider the space $X^\times$ as in Theorem \ref{fotetel1}   
or equivalently, in Theorem \ref{fotetel2}. Then there exists a smooth
Lorentzian metric $g_L$ on $X^\times$ such that $(X^\times, g_L)$ is, a
perhaps incomplete, Ricci-flat Lorentzian $4$-manifold.
\label{lorentztetel}
\end{theorem}

\noindent To make a comparison, let us indicate the ``size'' of the set of 
non-isometric solutions to the Lorentzian vacuum Einstein equation 
provided by Theorem \ref{lorentztetel}. 
By the fundamental classification result of Freedman \cite{fre}, 
connected and simply connected, oriented, closed topological four-manifolds 
are topologically classified by their intersection form 
$Q_M: H_2(M;\Z)\times H_2(M;\Z)\rightarrow H_4(M;\Z)\cong\Z$. By 
assumptions in our theorems here, $M$ is spin and smooth hence $Q_M$ must be 
even hence indefinite taking into account the other fundamental result in 
this field by Donaldson \cite{don}. Therefore if $\sigma (M)$ 
denotes the signature and $b_2(M)$ the second Betti number of $M$ then its 
intersection form looks like
\[Q_M=\frac{1}{8}\sigma (M)\left(\begin{smallmatrix}
2 & 0 &-1 & 0 & 0 & 0 & 0 & 0\\
0 & 2 & 0 &-1 & 0 & 0 & 0 & 0\\
-1& 0 & 2 &-1 & 0 & 0 & 0 & 0\\
0 &-1 &-1 & 2 &-1 & 0 & 0 & 0\\
0 & 0 & 0 &-1 & 2 &-1 & 0 & 0\\
0 & 0 & 0 & 0 &-1 & 2 &-1 & 0\\
0 & 0 & 0 & 0 & 0 &-1 & 2 &-1\\
0 & 0 & 0 & 0 & 0 & 0 &-1 & 2
\end{smallmatrix}\right)\oplus
\frac{1}{2}\left(b_2(M)-\sigma (M)\right)
\left(\begin{smallmatrix}0 &1\\1&0\end{smallmatrix}\right)\]
hence the simplest examples for $M$ are $S^4$, $S^2\times S^2$, the $K3$ 
surfaces, etc. Consequently, unlike the initial value formulation in the 
simply connected case, the set of solutions provided by Theorem 
\ref{lorentztetel} already contains many topologically different underlying 
spaces. But even more, most of these compact $M$'s 
themselves carry countable infinitely many different smooth structures, too. 
Finally, passing to the non-compact punctured spaces $M^\times$, the 
cardinality of the inequivalent smooth structures $X^\times$ underlying 
the Ricci-flat solutions in Theorem \ref{fotetel1} already reaches that of the 
continuum in ZFC set theory by a theorem of Gompf \cite{gom3} (recalled as 
Theorem \ref{gompftetel} here). Therefore the set of 
non-isometric Ricci-flat spaces exhibited in Theorem \ref{lorentztetel} is 
huge indeed. These solutions are not accessible within the initial value 
formulation because they, compared to the time evolution of typical initial 
data sets, are ``too long'' in an appropriate sense 
(cf. \cite[Section 5]{ete2}). Informally speaking, the vacuum Einstein 
equation is more tractable in Riemannian signature because of the elliptic 
nature of the Ricci-flatness condition in contrast to its 
hyperbolic character in Lorentzian signature: meanwhile solutions in 
Riemannian signature are protected by elliptic regularity hence ``extend 
well'', the regularity profiles of Lorentzian initial data sets quickly get 
destroyed during their hyperbolic time evolution. 

The paper is organized as follows. Section \ref{two} contains the 
collection of the required background material with rapid discussions of 
these results from our viewpoint. Sections \ref{three} and 
\ref{four}, respectively, contain the construction in the simpler non-exoting 
setting i.e., when $X^\times$ is not only homeomorphic but even diffeomorphic 
to $M^\times\subset M$ and then in the exotic setting with appropriate 
modifications. In Section \ref{five} we prove Theorem \ref{lorentztetel} 
by simply recalling \cite[Lemma 4.2]{ete2}. Finally 
in Section \ref{six} a physical interpretation of these Lorentzian Ricci-flat 
solutions is discussed. This interpretation places these solutions into 
the realm of the {\it strong cosmic censorship conjecture} and gravitational 
{\it topology change} processes.

%%%%%%%%%%%%%%%%%%%%%%%%%%%%%%%%%%%

\section{Background material}
\label{two}

%%%%%%%%%%%%%%%%%%%%%%%%%%%%%%

Let us begin with recalling all the powerful results, techniques, tools 
to be used during the construction of Riemannian Ricci-flat metrics 
in this paper.
 
{\it Construction of self-dual spaces}. 
It is well-known that the Fubini--Study metric on the complex projective 
space $\C P^2$ with orientation inherited from its complex structure is 
self-dual (or half-conformally flat) i.e., the anti-self-dual part $W^-$ of its 
Weyl tensor vanishes; consequently the oppositely oriented complex projective 
plane $(\C P^2)^{op}$ is anti-self-dual. A powerful generalization of 
this latter classical fact is Taubes' construction of an abundance of 
anti-self-dual $4$-manifolds; firstly we exhibit his result but now in an 
orientation-reversed form:

\begin{theorem} {\rm (Taubes \cite[Theorem 1.1]{tau2})} Let $M$ be a 
connected, compact, oriented smooth 
$4$-manifold. Let $\C P^2$ denote the complex projective plane with its 
usual orientation and let $\#$ denote the operation of taking the connected 
sum of manifolds. Then there exists a natural number $k_M\geqq 0$ such 
that for all $k\geqq k_M$ the modified compact manifold 
\[M\#\underbrace{\C P^2\#\dots\#\C P^2}_{k}\] 
admits a self-dual Riemannian metric. $\Diamond$

\label{taubestetel}
\end{theorem}

\noindent Let us roughly summarize how Taubes' construction works 
(\cite[Section 2]{tau2}). Take an arbitrary connected, 
oriented, closed Riemannian $4$-manifold $(M,g)$ and consider the 
density of the anti-self-dual part of the Weyl curvature of $g$ 
i.e., the pointwise norm $\vert W^-_g\vert_g$ along $M$. If it happens that 
somewhere around a point $p\in M$ this curvature density is large then take a 
$\C P^2$ with its usual Fubini--Study metric having zero anti-self-dual Weyl 
tensor and glue it to a ball $B^4_\varepsilon (p)\subset M$ of sufficiently 
small radius about the point. The result is a Riemannian metric on $M\#\C P^2$ 
having a bit smaller anti-self-dual Weyl tensor: this is because while 
$W^-_g$ is unchanged on $M\setminus \overline{B}^4_\varepsilon (p)$ it 
is killed in the bulk of $B^4_\varepsilon (p)$ except possibly along an annulus 
where $g$ and the Fubini--Study metric of $\C P^2$ have been glued together. 
Repeating this procedure, without doing connected summing on any previously 
added $\C P^2$ factor, probably very (but surely finitely) many times one 
comes up with a metric $\overline{\gamma}'$ on 
$\overline{X}_M:=M\#\C P^2\#\dots\#\C P^2$ (regarding the specific notation 
cf. Sections \ref{three} and \ref{four} below) whose 
$W^-_{\overline{\gamma}'}$ is already arbitrarily small in e.g. the original 
$L^2$-norm. Then, by the aid of the implicit 
function theorem, one perturbs this metric with a small symmetric tensor 
field $h$ on $\overline{X}_M$ into a new one 
$\overline{\gamma}:=\overline{\gamma}'+h$ which is already self-dual i.e. 
having $W^-_{\overline{\gamma}}=0$ along $\overline{X}_M$. For further 
rather technical details we refer to \cite[Section 2]{tau2}.

{\it Tools from twistor theory}. Let us now 
recall Penrose' twistor method \cite{pen1} to solve the Riemannian vacuum 
Einstein equation (for a very clear introduction cf. \cite[Chapter 13]{bes}, 
\cite{hit1,hit2}). Consider the bundle of unit-length anti-self-dual $2$-forms 
$S(\wedge^-\overline{X}_M)$ over a compact oriented space 
$(\overline{X}_M,\overline{\gamma})$ which is self-dual with respect to its 
orientation. Since in $4$ dimensions $\wedge^-\overline{X}_M$ 
is a rank $3$ real vector bundle over $\overline{X}_M$, its 
unit-sphere bundle $S(\wedge^-\overline{X}_M)$ is the total space of 
a smooth $S^2$-fibration 
$\overline{p}:S(\wedge^-\overline{X}_M)\rightarrow\overline{X}_M$. 
The Levi--Civita connection of the metric $\overline{\gamma}$ on 
$\overline{X}_M$ can be used to furnish the real $6$-manifold 
$S(\wedge^-\overline{X}_M)$ with a canonical almost complex structure; the 
fundamental observation of twistor theory is that this almost complex 
structure is integrable because $\overline{\gamma}$ is self-dual 
\cite[Theorem 13.46]{bes}. The resulting complex $3$-manifold 
$\overline{Z}\cong S(\wedge^-\overline{X}_M)$ is called the 
{\it twistor space} while the smooth fibration $\overline{p}:\overline{Z}
\rightarrow\overline{X}_M$ the {\it twistor fibration} of 
$(\overline{X}_M,\overline{\gamma})$. The most
important property of a twistor space of this kind is that its twistor
fibers $\overline{p}^{-1}(x)\subset\overline{Z}$ for all 
$x\in\overline{X}_M$ fit into
a locally complete complex $4$-paremeter family $\overline{X}_M^\C$ of
projective lines $Y\subset\overline{Z}$ each with normal bundle $H\oplus H$,
with $H$ being the dual of the tautological line bundle over $Y\cong\C P^1$. 
Moreover, there exists a real structure $\overline{\tau} 
:\overline{Z}\rightarrow \overline{Z}$ defined by 
taking the antipodal maps along the twistor fibers 
$Y_x:=\overline{p}^{-1}(x)\subset\overline{Z}$ for all 
$x\in\overline{X}_M\subset\overline{X}_M^\C$ which are 
therefore called {\it real lines} among all the lines $Y$ in $\overline{Z}$. 
In other words, $\overline{Z}$ is fibered exactly by the real lines $Y_x$ for 
all $x\in\overline{X}_M$. Hence the real $4$ dimensional self-dual geometry 
has been encoded into a $3$ dimensional complex analytic structure in the 
sense that one can recover $(\overline{X}_M,\overline{\gamma})$ just from 
$\overline{Z}$ up to conformal equivalence.

One can go further and raise the question how to recover precisely
$(\overline{X}_M,\overline{\gamma})$ itself from its conformal 
class, or more interestingly to us: how to get a Ricci-flat Riemannian 
$4$-manifold $(X_M,g)$ i.e., a solution of the (self-dual) 
Riemannian vacuum Einstein equation. Not surprisingly, to get the
latter stronger structure, one has to specify further data on the twistor
space. A fundamental result of twistor theory \cite{pen1} is that a 
solution of the $4$ dimensional (self-dual) Riemannian vacuum Einstein 
equation is equivalent to the following set of data (cf. 
\cite{hit1,hit2}):

\begin{itemize}

\item[$*$] A complex $3$-manifold $Z$, the total space of
a holomorphic fibration $\pi:Z\rightarrow\C P^1$;

\item[$*$] A complex $4$-paremeter family of holomorphically
embedded complex projective lines $Y\subset Z$, each with normal 
bundle $NY\cong H\oplus H$ (here $H$ is the dual of the tautological
bundle i.e., the unique holomorphic line bundle on $Y\cong\C P^1$ with
$\langle c_1(H), [Y]\rangle =1$);

\item[$*$] A non-vanishing holomorphic section $s$ of
$K_Z\otimes\pi^*H^4$ (here $K_Z$ is the canonical bundle of $Z$);

\item[$*$] A real structure $\tau :Z\rightarrow Z$ such that it 
coincides with the antipodal map $u\mapsto -\overline{u}^{-1}$ of $\C P^1$ 
upon restricting to the $\tau$-invariant elements $Y\subset Z$ 
(called real lines) from the family; moreover these real lines are both 
sections of $\pi$ and comprise a fibration of $Z$.

\end{itemize}

\noindent These data allow one to construct a Ricci-flat
and self-dual (i.e., the Ricci and the anti-self-dual Weyl part of
the curvature tensor vanishes) solution $(X_M,g)$ of the 
{\it Riemannian} Einstein's vacuum equation with vanishing cosmological 
constant as follows. The holomorphic lines $Y\subset Z$ form a locally 
complete family and fit together into a complex $4$-manifold 
$X_M^\C$. This space carries a natural complex conformal structure 
by declaring two nearby points $y_1,y_2\in X_M^\C$ to be 
null-separated if the corresponding lines intersect i.e., 
$Y_1\cap Y_2\not=\emptyset$ in $Z$. Infinitesimally this intersection 
condition means that on every tangent space $T_yX_M^\C\cong\C^4$ a null cone is
specified: using the identification $T_yX_M^\C\cong H^0(Y_y;\co(NY_y))
\cong H^0(\C P^1; \co (H\oplus H))$ given by $(a,b,c,d)\mapsto (au+b,cu+d)$, 
a tangent vector at $y$ is null if and only if its corresponding holomorphic 
sections have a common zero i.e. $ad-bc=0$ which is an equation of a cone. 
Restricting the complex conformal structure to the real lines singled out 
by $\tau$ and parameterized by an embedded real
$4$-manifold $X_M\subset X_M^\C$ we obtain the real 
conformal class $[g]$ of a Riemannian metric on $X_M$. 
The isomorphism $s:\pi^*H^{-4}\cong K_Z$ is essentially uniquely 
fixed by its compatibility with $\tau$ and gives rise to a volume form on 
$X_M$ this way fixing the metric $g$ in the conformal class. Given the 
conformal class, it is already meaningful to talk about the 
unit-sphere bundle of anti-self-dual $2$-forms $S(\wedge^-X_M)$ over $X_M$ 
with its induced orientation from the twistor space and $Z$ 
can be identified with the total space of $S(\wedge^-X_M)$. This way we 
obtain a smooth twistor fibration $p:Z\rightarrow X_M$ whose fibers are 
$\C P^1$'s hence $\pi :Z\rightarrow\C P^1$ can be regarded as a parallel 
translation along this bundle over $X_M$ with respect to a 
flat connection which is nothing but the induced connection of $g$ on 
$\wedge^-X_M$, cf. \cite{leb}. Knowing the decomposition of the Riemannian 
curvature into irreducible components over an oriented Riemannian $4$-manifold 
\cite{sin-tho}, this partial flatness of $S(\wedge^-X_M)$ implies that $g$ is 
Ricci-flat and self-dual. Finally note that, compared with the bare twistor 
space $\overline{Z}$ of a self-dual manifold 
$(\overline{X}_M,\overline{\gamma})$ above, the essential new 
requirement for constructing a self-dual {\it Ricci-flat} space $(X_M,g)$ 
is the existence of a holomorphic map $\pi$ from the 
twistor space $Z$ into $\C P^1$ which is compatible with the real 
structure in the above sense. We conclude our summary of the 
non-linear graviton construction by referring to \cite{hit1, hit2, hug-tod, 
leb, war-wel} or \cite[Chapter 13]{bes} for further details.

{\it Removable singularities in Yang--Mills fields}. Next let 
us refresh Uhlenbeck's by-now classical singularity removal theorem: 

\begin{theorem}{\rm (Uhlenbeck \cite[Theorem 4.1]{uhl} or 
\cite[Appendix D]{fre-uhl})} 

$*$\:Local version: Let $\nabla^\times$ be a solution of the 
${\rm SU}(2)$ Yang--Mills equations in the open punctured $4$-ball 
$B^4\setminus\{0\}$ with $\Vert F_{\nabla^\times}\Vert^2_{L^2(B^4)}
=\int_{B^4}\vert F_{\nabla^\times}\vert^2<+\infty$ i.e., having finite energy 
and $\nabla^\times =\dd +A^\times$ such that 
$A^\times\in L^2_1(B^4\setminus\{0\})$. Then $\nabla^\times$ is $L^2_2$ 
gauge equivalent to a connection $\nabla$ which extends smoothly across the 
singularity to a smooth connection.  

$*$\: Global version: Let $(M,g)$ be a connected, closed, oriented 
Riemannian $4$-manifold and let $\nabla^\times$ be an 
${\rm SU}(2)$ connection on a vector bundle $E^\times$ over 
$M^\times:=M\setminus\{{\rm point}\}$ which is a solution of the 
${\rm SU}(2)$ Yang--Mills equations and satisfies $\Vert F_{\nabla^\times}
\Vert_{L^2 (M)}<+\infty$ and there is an $L^2_{1,loc}$ gauge for 
$\nabla^\times$ around the puncturing of $M$. Then $\nabla^\times$ is 
$L^2_{2,loc}$ gauge equivalent to a connection $\nabla$ on a vector 
bundle $E$ over $M$ i.e., to a connection which extends across the pointlike 
singularity of the original connection. $\Diamond$  
\label{uhlenbecktetel}
\end{theorem}

\noindent Locally finite energy i.e., $F_{\nabla^\times}\in L^2_{loc}$ 
does not guarantee the continuity of the gauge transformation hence the 
topology of $E^\times$ can change i.e., $E^\times$ and $E$ can be different; 
however if $F_{\nabla^\times}\in L^{2+\varepsilon}_{loc}$ holds then we can 
assume continuity. Nevertheless the isomorphism class of $E$ is fully 
determined by the smooth connection $\nabla$ via the numerical value of the 
integral 
$-\infty <\frac{1}{8\pi^2}\int_M{\rm tr} (F_\nabla\wedge F_\nabla )<+\infty$, 
the second Chern number of the bundle $E$. 

{\it Exotic stuff}. Finally we evoke some results which provide 
us with a sort of summary of what is so special in four dimensions (i.e., 
absent in any other ones). First we recall a special class of large exotic 
(or fake) $\R^4$'s whose properties we will need here are summarized 
as follows:

\begin{theorem} {\rm (Gompf--Taubes, cf. \cite[Lemma 9.4.2, Addendum 9.4.4
and Theorem 9.4.10]{gom-sti})} There exists a pair $(R^4,K)$ consisting of a
differentiable $4$-manifold $R^4$ homeomorphic but not diffeomorphic to the
standard $\R^4$ and a compact oriented smooth $4$-manifold $K\subset R^4$
such that
\begin{itemize}

\item[$*$] $R^4$ cannot be smoothly embedded into the standard
$\R^4$ i.e., $R^4\not\subseteqq\R^4$ but it can be smoothly embedded as a
proper open subset into the complex projective plane i.e.,
$R^4\subsetneqq\C P^2$;

\item[$*$] Take a homeomorphism $f:\R^4\rightarrow R^4$, let
$0\in B^4_t\subset\R^4$ be the standard open $4$-ball of radius $t\in\R_+$
centered at the origin and put $R^4_t:=f(B^4_t)$ and
$R^4_{+\infty}:=R^4$. Then
\[\left\{ R^4_t\:\left\vert\:\mbox{$r\leqq t\leqq +\infty$ such that
$0<r<+\infty$ satisfies $K\subset R^4_r$}\right.\right\}\]
is an uncountable family of nondiffeomorphic exotic $\R^4$'s none of them
admitting a smooth embedding into $\R^4$ i.e., $R^4_t\not\subseteqq\R^4$
for all $r\leqq t\leqq +\infty$. $\Diamond$

\end{itemize}
\label{egzotikusnagycsalad}
\end{theorem}

\noindent 
The fact that any member $R^4_t$ in this family is not diffeomorphic to 
$\R^4$ implies the counterintuitive phenomenon that 
$R^4_t\not\cong W\times\R$ i.e., $R^4_t$ does not admit any {\it smooth} 
splitting into a $3$-manifold $W$ and $\R$ (with their unique smooth
structures) in spite of the fact that such {\it continuous} splittings
obviously exist. Indeed, from the contractibility of $R^4_t$ we can see that
$W$ must be a contractible open $3$-manifold (a so-called
{\it Whitehead continuum} \cite{whi}) however, by an early result of 
McMillen \cite{mcm} spaces of this kind always satisfy $W\times\R\cong\R^4$ 
i.e., their product with a line is always diffeomorphic to the standard 
$\R^4$. We will call this property of (any) exotic $\R^4$ occasionally below 
as ``creased''. 

From Theorem \ref{egzotikusnagycsalad} we deduce that for all 
$r<t<+\infty$ there is a sequence of smooth proper embeddings 
\[R^4_r\subsetneqq R^4_t\subsetneqq R^4_{+\infty}=R^4\subsetneqq\C P^2\] 
which are very wild in the following sense. The complement $\C P^2 
\setminus R^4$ of the largest member $R^4$ of this family is homeomorphic 
to $S^2$ regarded as an only ``continuously embedded projective line'' in 
$\C P^2$; therefore we shall denote this complement as $S^2:=\C 
P^2\setminus R^4\subset\C P^2$ in order to distinguish it from the 
ordinary projective lines $\C P^1=\C P^2\setminus\R^4 \subset\C P^2$. If 
$\C P^2=\R^4\cup\C P^1=\C^2\cup\C P^1$ is any holomorphic decomposition 
then $R^4\cap\C P^1\not=\emptyset$ (because otherwise $R^4\subseteqq\R^4$ 
would hold, a contradiction) as well as $S^2\cap \C P^1\not=\emptyset$ 
(because otherwise $H_2(R^4;\Z )\cong\Z$ would hold since $\C P^1\subset\C 
P^2$ represents a generator of $H_2(\C P^2;\Z)\cong\Z$, a contradiction 
again). Hence an ordinary projective line $\C P^1$ is always intersected by 
both $R^4$ and $S^2$ such that $S^2\cap\C P^1$ in the worst situation is a 
Cantor set. These demonstrate that the members of the large radial family 
``live somewhere between'' $\R^4$ and its complex projective closure $\C P^2$. 
However a more precise identification or location of them is a difficult task 
because these large exotic $\R^4$'s---although being honest differentiable 
$4$-manifolds---are very transcendental objects, cf. \cite[p. 366]{gom-sti}: 
They require infinitely many $3$-handles in any handle decomposition (like any 
other known large exotic $\R^4$) and there is presently\footnote{More precisely 
in the year 1999, cf. \cite{gom-sti}.} no clue as how one might draw 
explicit handle diagrams of them (even after removing their $3$-handles).
 
We note that the structure of small exotic $\R^4$'s i.e., which admit 
smooth embeddings into $\R^4$, is better understood and is quite different, cf. 
\cite[Chapter 9]{gom-sti}. For instance, unlike the large case, 
in their corresponding radial family certain (but surely not more than 
countably many) members are diffeomorphic such that the non-diffeomorphic 
small exotic $\R^4$'s are parameterized not by an interval but a Cantor 
set only, cf. \cite[Theorem 9.4.12 and its proof]{gom-sti}.  

Our last ingredient is the following {\it m\'enagerie} result of Gompf.

\begin{theorem} {\rm (Gompf \cite[Theorem 2.1]{gom3})} Let $X$ be a
connected (possibly non-compact, possibly with boundary) topological
$4$-manifold and let $X^\times:=X\setminus\{{\rm point}\}$ be
the punctured manifold with a single point removed. Then the non-compact
space $X^\times$ admits noncountably many (with the cardinality of the
continuum in ZFC set theory) pairwise non-diffeomorphic smooth
structures. $\Diamond$

\label{gompftetel}
\end{theorem}

\noindent If for instance $M$ is a connected compact smooth
$4$-manifold then Gompf's construction simply goes as follows: Take $R^4$ 
from Theorem \ref{egzotikusnagycsalad} and put 
\[X^\times:=M\#R^4\] 
which is a smooth $4$-manifold obviously homeomorphic 
to the punctured $M^\times$. More generally, the 
construction $X^\times_t:=M^\times\# R^4_t$ produces uncountably many 
mutually non-diffeomorphic smooth structures on 
the unique topological $4$-manifold underlying $X^\times_t$. 

%%%%%%%%%%%%%%%%%%%%%%%%%%%%%%%%%

\section{The construction}
\label{three}

%%%%%%%%%%%%%%%%%%%%%%%%%%%%%%%%%%%%%%%%%%%

\noindent In this section, which serves as a warming-up for the next 
one, we construct solutions of the vacuum Einstein equation on punctured 
$4$-manifolds carrying their standard smooth structure. We begin with 
an application of Theorem \ref{taubestetel} as follows. 

\begin{lemma}

Out of any connected, closed (i.e., compact without boundary) oriented 
smooth $4$-manifold $M$ one can construct a connected, open (i.e., 
non-compact without boundary) oriented smooth Riemannian $4$-manifold 
$(X_M,\gamma)$ which is self-dual but incomplete in general. 
\label{kezdolemma}
\end{lemma}

\noindent {\it Proof.} Pick any connected, oriented, closed, smooth 
$4$-manifold $M$. Referring to Theorem \ref{taubestetel} let 
$k:=\max (1,k_M)\in\N$ be a positive integer, put  
\[\overline{X}_M:=M\#\underbrace{\C P^2\#\dots\#\C P^2}_{k}\]
and let $\overline{\gamma}$ be a self-dual metric on it. Then 
$(\overline{X}_M,\overline{\gamma})$ is a compact self-dual 
manifold. Pick a $\C P^2$ factor within $\overline{X}_M$ and any 
(holomorphically embedded) projective line $\C P^1\subset\C P^2$ in that 
factor (avoiding its attaching point to $M$); then 
$\C P^1=\C P^2\setminus\C^2\cong\C P^2\setminus\R^4$ i.e., the line arises as 
the complement of an $\R^4$ in $\C P^2$. Let $K\subset\R^4$ be any connected 
compact subset and put 
\begin{equation}
X_M:=M\#\underbrace{(\C P^2\setminus\C P^1)\#\dots\#(\C 
P^2\setminus\C P^1)}_{k-1}\#_K(\C P^2\setminus\C P^1)
\cong M\#\underbrace{\R^4\#\dots\#\R^4}_{k-1}\#_K\R^4\cong 
M^\times\#\underbrace{\R^4\#\dots\#\R^4}_{k-1}
\label{ujsokasag}
\end{equation}
where the operation $\#_K$ means that the attaching point $y_0\in\R^4$ 
taken to glue a distinguished $\R^4$ with the rest 
$M\#\underbrace{\R^4\#\dots\#\R^4}_{k-1}$ 
satisfies $y_0\in K\subset\R^4$ and $M^\times:=
M\#_K\R^4\cong M\setminus\{{\rm point}\}$ is the punctured space with its 
inherited smooth structure from the smooth embedding $M^\times\subset M$. 
The result is a connected, open $4$-manifold $X_M$ (see Figure 1). 
\vspace{1in}

\centerline{
\begin{tikzpicture}[scale=0.7]
\node at (-0.8,0) {$M$};
\draw [thick] plot [tension=0.8,smooth cycle] coordinates {(0,0) 
(2,2) (5,0) (2,-2)};
\draw [thick] plot [tension=0.8,smooth] coordinates {(0.5,0.4) (1.5,0.2) 
(2.5,0.4)};
\draw [thick] plot [tension=0.8,smooth] coordinates {(0.8,0.3) (1.5,0.5) 
(2.2,0.3)};
\node at (7.2,0) {$X_M$};
\draw [thick] plot [tension=0.8,smooth] coordinates { (8.5,0.2) (9,-0.2)
(9.5,0.2)};
\draw [thick] plot [tension=0.8,smooth] coordinates { (8.7,0) (9,0.1) 
(9.3,0)};
\draw [thick] plot [tension=0.8,smooth] coordinates
{(8.5,1) (9,1.3) (10,1) (11,0.7) (12,0.9) (12.5,1.3)};
\draw [thick] plot [tension=0.8,smooth] coordinates
{(13.1,1.3) (13.4,0.6) (14,0.3)};
\draw plot [tension=0.8,smooth] coordinates
{(12.5,1.3) (12.8,1.4) (13.1,1.3)};
\draw plot [tension=0.8,smooth] coordinates
{(12.5,1.3) (12.8,1.2) (13.1,1.3)};
\draw [thick] plot [tension=0.8,smooth] coordinates
{(8.5,-1) (9,-1.3) (10,-1) (11,-0.7) (12,-0.9) (12.5,-1.3)};
\draw [thick] plot [tension=0.8,smooth] coordinates
{(13.1,-1.3) (13.4,-0.6) (14,-0.3)};
\draw plot [tension=0.8,smooth] coordinates
{(12.5,-1.3) (12.8,-1.4) (13.1,-1.3)};
\draw plot [tension=0.8,smooth] coordinates
{(12.5,-1.3) (12.8,-1.2) (13.1,-1.3)};
\draw [thick] plot [tension=0.8,smooth] coordinates {(8.5,1) (8,0)
(8.5,-1)};
\draw [thick] plot [tension=0.8,smooth] coordinates {(14,0.3) (15,0.2) 
(16,1)};
\draw [thick] plot [tension=0.8,smooth] coordinates {(14,-0.3) (15,-0.2) 
(16,-1)};
\draw [gray] plot [tension=0.8,smooth] coordinates {(16,1) (15.8,0) (16,-1)};
\draw [gray] plot [tension=0.8,smooth] coordinates {(16,1) (16.2,0) (16,-1)};
%\label{abra1}
\end{tikzpicture}}

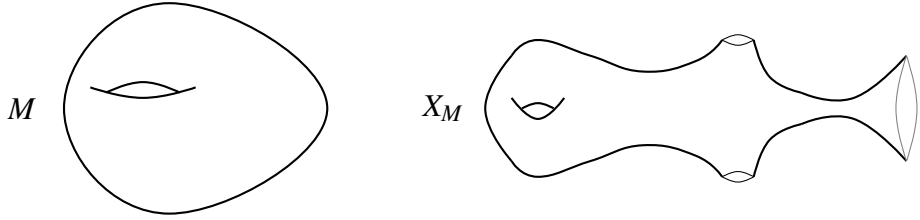
\centerline{Figure 1. Construction of $X_M$ out of $M$. The gray 
ellipse represents a distinguished end} 
\centerline{diffeomorphic to the complement of a connected compact subset  
$K$ in $\R^4$.} 
\vspace{0.3in}

\noindent From the proper smooth embedding $X_M\subsetneqq\overline{X}_M$ 
there exists a restricted self-dual Riemannian metric
$\gamma:=\overline{\gamma}\vert_{X_M}$ on $X_M$ which is however in
general non-complete. $\Diamond$
\vspace{0.1in}

\noindent Next we improve the incomplete 
self-dual space $(X_M,\gamma)$ of Lemma \ref{kezdolemma} to a complete 
Ricci-flat space $(X_M,g)$ by conformally rescaling $\gamma$ 
with a suitable positive smooth function $\varphi :X_M\rightarrow\R_+$ which 
is a ``multi-task'' function in the sense that it kills both the scalar 
curvature and the traceless Ricci tensor of $\gamma$ moreover blows up 
sufficiently fast along the $\R^4$ ends of $X_M$ to render the 
rescaled metric $g$ complete. Two classical examples serve as a
motivation.
\vspace{0.1in}

\noindent {\it First example}. 
First, let $S^4\subset\R^5$ be the standard $4$-sphere equipped with the 
standard orientation and round metric inherited from the embedding. 
Put $\overline{X}_M:=S^4$ and $\overline{\gamma}:=$the standard 
round metric. It is well-known that $(\overline{X}_M,\overline{\gamma})=
(S^4,\overline{\gamma})$ is self-dual and Einstein with non-zero 
cosmological constant i.e., not Ricci-flat. Put 
$X_M:=S^4\setminus\{\infty\}=\R^4$; then 
$\gamma=\overline{\gamma}\vert_{\R^4}$ thus $(X_M,\gamma )=(\R^4,\gamma )$ 
is an incomplete self-dual space. But setting $\varphi 
:\R^4\rightarrow\R_+$ to be $\varphi (x):=(1+\vert x\vert^2)^{-1}$, 
then $g:=\varphi^{-2}\gamma$ 
is nothing but the standard flat metric $\eta$ on $\R^4$ which is of 
course complete and Ricci-flat. Hence $(X_M,g)=(\R^4, \eta)$, 
the conformal rescaling of $(X_M,\gamma )=(\R^4,\gamma)$, is the 
desired complete Ricci-flat space in this simple case. Note that 
$(\R^4,\eta )$ is a trivial hyper-K\"ahler space, too.

It is worth working out here how the corresponding holomorphic map 
$\pi :Z\rightarrow \C P^1$ over the corresponding twistor space 
arises in this situation (see the summary of twistor theory in Section 
\ref{two}). Consider the smooth twistor fibration 
$\overline{p}:\overline{Z}\rightarrow S^4$. Since $\R^4\subset S^4$ writing 
$Z:=\overline{Z}\vert_{\R^4}$ and $p:=\overline{p}\vert_{\R^4}$ we obtain a 
restricted fibration $p:Z\rightarrow\R^4$. Unlike the full twistor fibration 
over $S^4$, the restricted one is topologically trivial i.e. $Z$ is 
homeomorphic to $\R^4\times S^2$ since $\R^4$ is contractible; consequently 
$Z$ admits a continuous trivialization over $\R^4$. This is a necessary
topological condition for the existence of the map $\pi$. 
Since $\R^4$ with its flat metric is conformally 
equivalent to $S^4\setminus\{\infty\}$ with its round metric, $Z$ 
arises by deleting the twistor line over $\infty\in S^4$ from $\overline{Z}$. 
However it is well-known that the twistor space $\overline{Z}$ of the round 
$S^4$ is $\C P^3$ consequently the twistor space $Z$ of the flat $\R^4$ is 
simply $\C P^3\setminus\C P^1$. More explicitly, 
take a homogeneous coordinate system $[z_0:z_1:z_2:z_3]$ on $\C P^3$ and 
remove the line $z_0=z_1=0$ from $\C P^3$ to get $Z$. We wish to define 
a map $\pi :Z\rightarrow\C P^1$ such that its target space is a twistor 
i.e. a real line in $Z$. Any line in $\C P^3\setminus\C P^1$ can be 
written as $[z_0:z_1:az_1+bz_0:cz_1+dz_0]$ with $[z_0:z_1]\in\C P^1$ and 
$a,b,c,d\in\C$ being some parameters. Note that the case of $a=b=c=d=0$ is 
meaningful and $[z_0:z_1:0:0]$ is simply the distinguished line 
$[z_0:z_1]$ in $\C P^3\setminus\C P^1$. Thus
\begin{equation}
\mbox{The lines in $Z$}=\{[z_0:z_1:az_1+bz_0:cz_1+dz_0]\:\vert\: 
\mbox{$[z_0:z_1]\in\C P^1$ and $a,b,c,d\in\C$}\}\:\:.
\label{H+H}
\end{equation}
The real structure on $\overline{Z}$ is defined by demanding the fibers 
of $\overline{p}:\overline{Z}\rightarrow S^4$ to be invariant. Under 
$\overline{Z}\cong\C P^3$ it comes from the identification $\C^4\cong\HH^2$ 
and has the form $[z_0:z_1:z_2:z_3]\mapsto 
[\overline{z}_1:-\overline{z}_0:\overline{z}_3:-\overline{z}_2]$. It is 
compatible with the antipodal map $[z_0:z_1]\mapsto 
[\overline{z}_1:-\overline{z}_0]$ and restricts to a real structure 
$\tau:Z\rightarrow Z$. It then follows that the corresponding real lines 
have the shape $[z_0:z_1:az_1+\overline{c}z_0:cz_1-\overline{a}z_0]$. 
Consequently the twistor fibration $p:Z\rightarrow\R^4$ looks like
$[z_0:z_1:az_1+\overline{c}z_0:cz_1-\overline{a}z_0]\mapsto (a,c)
\in\C^2\cong\R^4$ and in particular the distinguished line is real and can 
be identified with the twistor line $p^{-1}(0)$ over the origin. Since every 
point $z\in Z$ contained in exactly one real line let us define $\pi  
:Z\rightarrow p^{-1}(0)$ by the canonical projection $\pi 
([z_0:z_1:az_1+\overline{c}z_0:cz_1-\overline{a}z_0]):=[z_0:z_1]$. 
Upon introducing the projective coordinate 
$u:=\frac{cz_1-\overline{a}z_0}{az_1+\overline{c}z_0}$ if $(a,c)\not=(0,0)$ or 
$u:=\frac{z_1}{z_0}$ if $(a,c)=(0,0)$ along the twistor lines in the 
domain of $\pi$ the map looks like 
\begin{equation}
\pi (u)=\left\{\begin{array}{ll}\frac{\overline{c}u+\overline{a}}{-au+c} 
& \mbox{if $(a,c)\not=(0,0)$}\\         
u & \mbox{if $(a,c)=(0,0)$}
              \end{array}\right.
\label{alappi}
\end{equation}
which is an obviously holomorphic map since it arises by holomorphic
deformations of $p^{-1}(0)$ within $Z$ moreover it is the identity on
$p^{-1}(0)$. What we only have to check is that $\pi$ is compatible with
the real structure. This means that we have to demonstrate that all real 
lines $p^{-1}(x)\subset Z$ are sections of $\pi :Z\rightarrow p^{-1}(0)$ 
or in other words that $\pi\vert_{p^{-1}(x)}: p^{-1}(x)\rightarrow p^{-1}(0)$ 
is a holomorphic bijection of $\C P^1$ for every $x\in\R^4$. Assume that
this is not true. Since $\pi\vert_{p^{-1}(x)}$ has the form (\ref{alappi}) 
we can normalize its coefficients such that $\vert c\vert^2+\vert 
a\vert^2=1$. However the assumption implies that this rational function is 
constant in $u$ yielding $\vert c\vert^2+\vert a\vert^2=0$, a contradiction.
\vspace{0.1in}

\noindent {\it Second example}. This time put $\overline{X}_M:=\C P^2$ and 
$\overline{\gamma}:=$Fubini--Study metric. It is well-known that 
$(\overline{X}_M,\overline{\gamma})=(\C P^2, 
\overline{\gamma})$ is self-dual and Einstein with non-zero cosmological 
constant i.e., not Ricci-flat. Now let $X_M:=\C P^2\setminus\C P^1=\R^4$; 
then $\gamma=\overline{\gamma}\vert_{\R^4}$ and $(X_M,\gamma )=(\R^4,\gamma )$ 
is an incomplete self-dual space. 
If $0\not=(z_0,z_1,z_2)\in\C^3$ and $[z_0:z_1:z_2]\in\C P^2$ then take 
the projective line $\C P^1\subset\C P^2$ defined by $z_0=0$. Introducing 
$w_i=\frac{z_i}{z_0}$ ($i=1,2)$ and $w=(w_1,w_2)\in\C^2=\R^4$, on the 
complementum $\C P^2\setminus\C P^1$ the restricted Fubini--Study metric 
$\gamma$ looks like 
\[\gamma_{ij}(w)=(1+\vert w\vert^2)^{-1}\delta_{ij}-
(1+\vert w\vert^2)^{-2}\overline{w}_iw_j\] 
and along this local part it already possesses a K\"ahler potential 
$K(w)=\log (1+\vert w\vert^2)$. 
This time define $\varphi :\R^4\rightarrow \R_+$ as $\varphi (w):=
{\rm e}^{-\frac{3}{4}K(w)}=(1+\vert w\vert^2)^{-\frac{3}{4}}$ which is a 
non-holomorphic function and consider the conformally rescaled (real) metric 
$g:=\varphi^{-2}\gamma$. One can check that this is a complete 
Ricci-flat metric on $\R^4$. Hence $(X_M,g)=(\R^4,g)$, the conformal 
rescaling of $(X_M,\gamma )=(\R^4,\gamma )$ is a complete Ricci-flat space. 
It is already not flat but note again that nevertheless $g$ indudes a (not 
asymptotically flat in any sense) hyper-K\"ahler structure on $\R^4$ 
because $g$ is a complete, self-dual, Ricci-flat metric on the simply 
connected space $\R^4$. 

Again, the corresponding twistor-theoretic map $\pi$ arises as follows. 
Consider the smooth twistor fibration $\overline{p}:\overline{Z}
\rightarrow\C P^2$. Since $\R^4\subset\C P^2$, writing 
$Z:=\overline{Z}\vert_{\R^4}$ and $p:=\overline{p}\vert_{\R^4}$
we obtain a restricted fibration $p:Z\rightarrow\R^4$.
Unlike the full twistor fibration over $\C P^2$, this restricted one is
topologically trivial i.e., $Z$ is homeomorphic to $\R^4\times S^2$
since $\R^4$ is contractible; consequently $Z$ admits a continuous
trivialization over $\R^4$. This is a necessary topological condition for
the existence of the map $\pi$. It is known that $\overline{Z}
\cong P(T\C P^2)$ i.e., the twistor space of the complex projective space can 
be identified with its projective holomorphic tangent bundle. Consequently
$\overline{Z}$ admits a very classical description namely
can be identified with the flag manifold $F_{12}(\C ^3)$ consisting of pairs
$(L,P)$ where $0\in L\subset\C^3$ is a
line (i.e., a point $p\in\C P^2$) and $0\in L\subset P\subset\C^3$ is a 
plane containing the line (i.e., a line $p\in\ell\subset\C P^2$ containing 
the point). Then in the twistor fibration 
$\overline{p}:\overline{Z}\rightarrow\C P^2$ of the complex projective space 
$(L,P)\in F_{12}(\C ^3)$ is sent into the point $x\in\C P^2$ 
provided by the line $X:=L^\perp\cap P\subset\C ^3$ where $L^\perp$ is the 
plane perpendicular to the line $L$ in $\C^3$ with respect to the standard 
Hermitian scalar product. This is a smooth but not holomorphic fibration 
over $\C P^2$ with $\C P^1$'s as fibers since $\overline{p}^{-1}(x)=
\{ (L,P)\:\vert\:L\subset X^\perp\:,\:X\subset P\}=\{(p,\ell )\:\vert\: 
x^\perp\cap\ell, x\in\ell\}$ i.e., it consists of 
all lines $\ell\subset\C P^2$ through $x\in\C P^2$ (a copy of $\C P^1$) 
and a distinguished point $p$ on each given by its intersection with the line 
$x^\perp\subset\C P^2$ given by $X^\perp\subset\C^3$. 
Consider now the restricted twistor fibration $p:Z\rightarrow\R^4$. Fix a 
point $x_0\in\C P^2 \setminus x_0^\perp =\R^4$ with target space 
$p^{-1}(x_0)\cong\C P^1$ consisting of terminating pairs
$(p_0,\ell_0)\in p^{-1}(x_0)\subset Z$. Take a starting pair $(p,\ell)\in Z$ 
over a running point $x\in \C P^2\setminus x^\perp_0$. 
Our aim is to construct a holomorphic map which associates to
$(p,\ell)$ another pair $(p_0,\ell_0)$. We construct this
$\pi:Z\rightarrow p^{-1}(x_0)$ very simply as follows.
Consider a starting pair $(p,\ell)$ and take its
line component $\ell\subset\C P^2$. This line has a unique
intersecion $p_0:=x_0^\perp\cap\ell$ with the infinitely distant 
projective line. Then, given the target space $p^{-1}(x_0)$, define the 
projective line component $\ell_0\subset\C P^2$ in the 
terminating pair $(p_0,\ell_0)\in p^{-1}(x_0)$ by taking
the unique projective line $\ell_0$ connecting
$p_0$ with $x_0$. In short,
\begin{eqnarray}
\mbox{$\pi((p,\ell)):=
(p_0,\ell_0)$ where}& &
\mbox{$p_0\in\C P^2$ satisfies $p_0:=x^\perp_0\cap\ell$ and}\nonumber\\
& &\mbox{$\ell_0\subset\C P^2$ satisfies that
$\ell_0$ connects $p_0$ with $x_0$ in $\C P^2$}
\label{pi}
\end{eqnarray}
(see Figure 2 for a construction of this map in projective geometry).

\centerline{
\setlength{\unitlength}{1cm}
\begin{picture}(7,5)
\thicklines
\put(1,1){\circle*{0.2}}
\put(1.2,0.8){$x$}
\put(3,3){\circle*{0.2}}
\put(2.8,2.5){$p_0$}
\put(5,1){\circle*{0.2}}
\put(4.4,0.8){$x_0$}
\put(-0.5,-0.5){\line(1,1){4}}
\put(1.3,1.7){$\ell$}
\put(5.5,0.5){\line(-1,1){3}}
\put(4.6,1.7){$\ell_0$}
\put(-1,3){\line(1,0){6}}
\put(1,3.3){$x^\perp_0$}
\put(0,-0.8){\line(0,1){5}}
\put(0,0){\circle*{0.2}}
\put(-0.4,0.1){$p$}
\put(-0.6,2){$x^\perp$}
\put(7,2){$\C P^2$}
%\label{abra2}
\end{picture}
}
\vspace{0.3in}

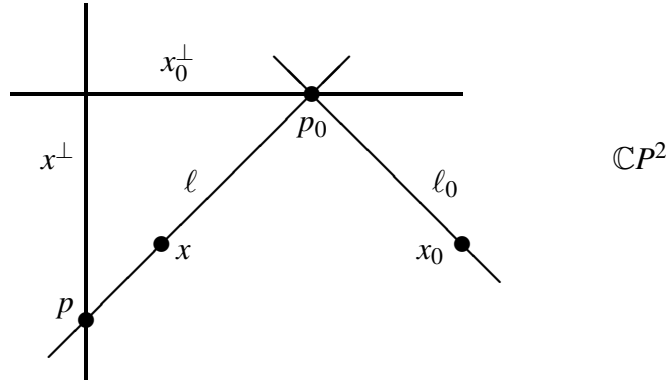
\centerline{Figure 2. Construction of the map $\pi : Z\rightarrow \C P^1$ 
satisfying $\pi((p,\ell))=(p_0,\ell_0)$.}
\vspace{0.3in}

\noindent It is a classical observation that this map is
well-defined on $Z$ and holomorphic; in particular it is the identity on the 
target space $p^{-1}(x_0)$ i.e., $\pi((p_0,\ell_0))=(p_0,\ell_0)$. 

This globally defined map admits a local description which looks very similar 
to the {\it First example}. Given the target point $x_0\in\C 
P^2\setminus x_0^\perp\cong\R^4$, its twistor line can be identified 
with the infinitely distant line $x^\perp_0\subset\C P^2$ or 
equivalently, $x^\perp_0\subset Z$. Likewise, if $x\in\C P^2\setminus 
x_0^\perp\cong\R^4$ is a nearby point then its twistor line is 
$x^\perp\subset\C P^2$ or equivalently, $x^\perp\subset Z$. In this 
picture the map (\ref{pi}) can be described simply as follows: 
if $x\in\ell\subset\C P^2$ is a line then $\pi (\ell\cap x^\perp )= 
\ell\cap x^\perp_0$ as $\ell$ runs over all possibilities. Pick 
homogeneous coordinates $[z_0:z_1:z_2]$ on $\C P^2$ such that 
$x_0:=[1:0:0]$ hence $x_0^\perp =\{[0:v_1:v_2]\vert [v_1:v_2]\in\C P^1\}$. 
Likewise, if $x=[1:z_1:z_2]$ is the nearby point then 
$x^\perp = \{[-\overline{z}_1w_1-\overline{z}_2w_2:w_1:w_2]\vert [w_1:w_2]\in\C 
P^1\}$. The affine part of the line $\ell$ connecting $[1:z_1:z_2]$ and 
$[0:v_1:v_2]$ is $\{[t:v_1+(z_1-v_1)t: v_2+(z_2-v_2)t]\:\vert\:t\in\C\}$ 
hence by solving the equation $[t:v_1+(z_1-v_1)t: v_2+(z_2-v_2)t]=
[-\overline{z}_1w_1-\overline{z}_2w_2:w_1:w_2]$ 
for $[w_1:w_2]$ and upon introducing the projective coordinate 
$u:=\frac{w^1}{w^2}$ along $x^\perp\cong \C P^1$ the map (\ref{pi}) takes 
the shape 
\[\pi (u)=\frac{(1+\vert z_1\vert^2)u-z_1\overline{z}_2}
{-\overline{z}_1z_2u+(1+\vert z_2\vert^2)}\]
hence looks like (\ref{alappi}) indeed.

\begin{remark}\rm It follows from the description (\ref{H+H}) of its 
holomorphic lines that the twistor space $Z$ of the flat $\R^4$ can be 
globally holomorphically identified with the total space of the bundle 
$H\oplus H$ over the distinguished projective line $\C P^1$ parameterized 
with $[z_0:z_1]$ in (\ref{H+H}) and the map (\ref{alappi}) 
is nothing but the projection $\pi :H\oplus H\rightarrow\C P^1$. The point is 
that this picture on the twistor space continues to hold true in the generic 
case at least locally. Consider a 
general twistor space $Z$ with its twistor fibration $p:Z\rightarrow X_M$. 
Take $x\in X_M$ and let $Np^{-1}(x)$ be the normal bundle of the twistor 
line $p^{-1}(x)\subset Z$. We know (see the summary of twistor theory in 
Section \ref{two}) that the holomorphy type of the normal bundle is fixed in 
advance and is a very special bundle: it is positive hence admits holomorphic 
sections such that they parameterize a locally complete family of projective 
lines $Y\subset Z$ which are small holomorphic 
deformations of $Y_x=p^{-1}(x)$ inside $Z$ (cf. e.g. \cite[Sections III.1 and 
III.2]{gri}) including therefore all nearby real lines as well. 
Thus there exist small open neighbourhoods $U_x\subset X_M$ of $x$ 
and $V_x\subset Np^{-1}(x)$ of the zero section with an injection 
$\Psi_x: p^{-1}(U_x)\rightarrow V_x$ that is, $\Psi_x$ maps 
injectively the twistor fibers over $U_x$ into the space of holomorphic 
sections of $V_x$ such that this map is onto an appropriately defined 
subspace of real sections. More explicitly, we know that the normal bundle is 
always isomorphic to $H\oplus H$ consequently all small holomorphic 
deformations of a given twistor line within $Z$ can be parameterized by 
$(a,b,c,d)\in\C^4\cong H^0(\C P^1;\co (H\oplus H))$ such that the twistor 
lines satisfy a reality condition implying $ad-bc\not=0$ (because the real 
lines never intersect), exactly like in the {\it First example}. 
\end{remark}

\noindent We can return now to the much more 
general situation set up in Lemma \ref{kezdolemma}; motiveted by 
the examples, instead of finding conformal rescalings 
$\varphi :X_M\rightarrow\R_+$ directly, we are going to use Penrose' 
non-linear graviton construction (i.e., twistor theory \cite{pen1}) 
to find their holomorphic counterparts $\pi :Z\rightarrow\C P^1$. Consider 
the compact self-dual space $(\overline{X}_M,\overline{\gamma})$ from Lemma 
\ref{kezdolemma}, take its twistor fibration $\overline{p}: 
\overline{Z}\rightarrow\overline{X}_M$ and let $p:Z\rightarrow X_M$ 
be its restriction induced by the smooth embedding $X_M\subsetneqq 
\overline{X}_M$ i.e., $Z:=\overline{Z}\vert_{X_M}$ and 
$p:=\overline{p}\vert_{X_M}$. Then $Z$ is a non-compact 
complex $3$-manifold already obviously possessing all the required twistor 
data except the existence of a holomorphic mapping $\pi :Z\rightarrow\C P^1$. 

\begin{lemma} Consider the connected, open, oriented, incomplete, 
self-dual space $(X_M,\gamma)$ as in Lemma \ref{kezdolemma} 
with its twistor fibration $p:Z\rightarrow X_M$ constructed above. 
If $\pi_1(M)=1$ and $M$ is spin (or equivalently, having even  
intersection form) then there exists a holomorphic mapping 
$\pi :Z\rightarrow\C P^1$.
\label{pilemma}
\end{lemma}

\noindent {\it Proof.} Let $x_0\in X_M$ be a fixed point. 
Our aim is to construct a holomorphic map 
\begin{equation}
\pi:\:\:Z\longrightarrow p^{-1}(x_0)\cong\C P^1
\label{celfuggveny}
\end{equation}
that we carry out by analytic continuation. 

First, put $\pi\vert_{p^{-1}(x_0)}:={\rm Id}_{p^{-1}(x_0)}$. Secondly, 
suppose that in $x\in X_M$ the map is already defined i.e. there exists 
$\pi\vert_{p^{-1}(x)}:  p^{-1}(x)\rightarrow p^{-1}(x_0)$ which is 
compatible with the real structure on $Z$ hence is a 
holomorphic bijection between the twistor fibers in question or in other 
words is a holomorphic bijection of $\C P^1$. Consider a sufficiently small 
open neighbourhood $U_x\subset X_M$ of $x$ such that $p^{-1}(U_x)\subset Z$ can 
be holomorphically modeled within the neighbourhood $V_x$ of the zero 
section of $Np^{-1}(x)$, the normal bundle of the twistor line 
$p^{-1}(x)$. Define $\rho_x: p^{-1}(U_x)\rightarrow p^{-1}(x)$ to be the 
restriction of the projection $\pi :Np^{-1}(x)\rightarrow p^{-1}(x)$ onto 
the image of the twistor lines of $p^{-1}(U_x)$ within $Np^{-1}(x)$. 
That is, given a point $z\in p^{-1}(U_x)$ there exists a unique real 
line passing through it and $\rho_x(z)\in p^{-1}(x)$ simply arises 
by the projection of this line onto the central twistor line $p^{-1}(x)$. 
This local map is clearly holomorphic because it stems from holomorphic 
deformations of $p^{-1}(x)$ inside $Z$ provided by its locally complete 
family of lines.\footnote{For a comparison with the general theory 
\cite[Proposition 1.3]{gri} we remark here that although the Griffiths 
obstruction groups $H^1( p^{-1}(x); \co ((\pi\vert_{p^{-1}(x)})^*Tp^{-1}(x_0)
\otimes S^kN^*p^{-1}(x)))$ against the extendibility of 
$\pi\vert_{p^{-1}(x)}:p^{-1}(x)\rightarrow p^{-1}(x_0)$ to the $k^{\rm th}$
formal neighbourhoud of $p^{-1}(x)\subset Z$ are non-trivial for $k\geqq 4$,
the above construction (or the explicit {\it Second example}) shows that 
the corresponding obstruction classes $\omega(\pi_{k-1})$ themselves are 
nevertheless trivial. This essentially follows from Kodaira's 
integrability condition $H^1(p^{-1}(x);\co (Np^{-1}(x)))=\{0\}$, cf. 
\cite[Theorem 3.1]{gri}.} Moreover $\rho_x$ is the identity on $p^{-1}(x)$. 
What we have to still check that it is compatible with the real 
structure on $Z$ i.e. for every $y\in U_x$ the map $\rho_x$ is a 
holomorphic bijection between $p^{-1}(y)$ and $p^{-1}(x)$. 
%Under the isomorphisms $H^0(p^{-1}(x);\co (Np^{-1}(x)))\cong 
%T_xX_M\otimes_\R\C\cong\Sigma_x^+\otimes\Sigma_x^-\cong{\rm Hom} 
%((\Sigma_x^-)^*;\Sigma_x^+)$ provided by the orientation, the metric $\gamma$ 
%and the spin structure on$X_M$, the twistor line $p^{-1}(y)$ with $y\in U_x$ 
%is mapped into a unitary map ${\rm U}_y: (\Sigma_x^-)^*\rightarrow \Sigma_x^+$ 
%between the chiral spinor bundles over $x\in X_M$. 
Exploiting the isomorphism $Np^{-1}(x)\cong H\oplus H$ (see the summary of 
twistor theory in Section \ref{two}) the map $\rho_x:  
p^{-1}(U_x)\rightarrow p^{-1}(x)$ can be described by the projection 
$\pi :H\oplus H\rightarrow\C P^1$ therefore, upon introducing the projective 
coordinate $u$ along $p^{-1}(y)\cong\C P^1$ 
\[\rho_x(u)=\left\{\begin{array}{ll}\frac{au+b}{cu+d}
& \mbox{if $(a,b,c,d)\not=(0,0,0,0)$}\\
         u & \mbox{if $(a,b,c,d)=(0,0,0,0)$}
              \end{array}\right.\]
where $(a,b,c,d)\in H^0(\C P^1;\co (H\oplus H))\cong\C^4$ are the 
coefficients of a real line hence satisfy an appropriate reality condition. 
Fortunately whatever this reality condition is, we surely know that 
$ad-bc\not=0$ because real lines never intersect. However this implies that 
the map $u\mapsto \frac{au+b}{cu+d}$ is not constant in $u$ that is, $\rho_x$ 
is indeed a holomorphic bijection between $p^{-1}(y)$ and $p^{-1}(x)$ 
for all $y\in U_x$ (such that it is the identity on $p^{-1}(x)$) hence 
$\rho_x$ is compatible with the real structure on $p^{-1}(U_x)$ as desired. 

Therefore let us define the local extension $\pi\vert_{p^{-1}(U_x)}:p^{-1}(U_x)
\rightarrow p^{-1}(x_0)$ by the composition 
$\pi\vert_{p^{-1}(U_x)}:=\pi\vert_{p^{-1}(x)}\circ\rho_x$. 
By assumption $\pi\vert_{p^{-1}(x)}$ already possesses all the 
required properties hence is compatible with the real structure therefore it 
is a holomorphic bijection between the twistor lines $p^{-1}(x)$ and 
$p^{-1}(x_0)$; consequently, taking a projective coordinate $v$ 
along $p^{-1}(x)\cong\C P^1$, we know that 
$\pi\vert_{p^{-1}(x)}$ also has the form 
$v\mapsto\frac{a_0v+b_0}{c_0v+d_0}$ with some 
$a_0,b_0,c_0,d_0\in\C$ satisfying $a_0d_0-b_0c_0\not=0$. Composing the maps 
above means that we insert $v=\frac{au+b}{cu+d}$ where $u$ is the projective 
coordinate along $p^{-1}(y)\cong\C P^1$ as before; thus the local extension 
looks like  
\[\pi\vert_{p^{-1}(U_x)}(u)=\frac{a_0\frac{au+b}{cu+d} +b_0}
{c_0\frac{au+b}{cu+d}+d_0}=\frac{(a_0a+b_0c)u+(a_0b+b_0d)}
{(c_0a+d_0c)u+(c_0b+d_0d)}\:\:.\]
Since $(a_0a+b_0c)(c_0b+d_0d)-(a_0b+b_0d)(c_0a+d_0c)=
(a_0d_0-b_0c_0)(ad-bc)\not=0$ it readily follows that it 
continues to be compatible with the real structure. 

Thirdly, since $M$ is connected, simply connected and spin, $Z$ is 
connected, simply connected and $p:Z\rightarrow X_M$ is trivial. These make 
sure that $\pi$ extends over $Z$ in a consistent way. $\Diamond$

\begin{remark}\rm Note that the reasons for both the local map $\rho_x$ 
and the non-local one $\pi\vert_{p^{-1}(x)}$ having the same shape (namely 
both are fractional linear transformations of $\C P^1$) are quite different. 
Nevertheless it makes possible to regard $\pi :Z\rightarrow\C P^1$ as an 
action of ${\rm SL}(2;\C )$ on the target projective line $\C P^1$ via 
fractional linear transformations which are in turn ${\rm SO}(3)$ 
rotations on $S^2$ regarded as the unit sphere in the space of 
anti-self-dual $2$-forms provided either by the old or the new metric 
$\gamma$ or $g$, respectively. 
\end{remark}

\noindent It follows that $\pi :Z\rightarrow\C P^1$ i.e., the map 
(\ref{celfuggveny}) constructed in Lemma \ref{pilemma} is compatible with the 
real structure $\tau :Z\rightarrow Z$ already fixed by the self-dual structure 
in Theorem \ref{taubestetel} therefore twistor theory provides us with a 
Ricci-flat (and self-dual) Riemannian metric $g$ on $X_M$. We proceed 
further and demonstrate that, unlike $(X_M,\gamma)$, the space $(X_M,g)$ is 
complete. 
\begin{lemma} The four dimensional connected and simply connected, open, 
oriented, Ricci-flat Riemannian spin manifold $(X_M,g)$ is complete. 
\label{teljeslemma}
\end{lemma}

\noindent {\it Proof}. Since both $\gamma$ and 
this Ricci-flat metric $g$ arise from the same complex structure on 
the same twistor space $Z$ we know from twistor theory that these metrics are 
in fact conformally equivalent. That is, there exists a smooth non-constant 
strictly positive function $\varphi :X_M\rightarrow\R_+$ such that 
$\varphi^{-2}\gamma= g$. Our strategy to prove completeness 
is to follow Gordon \cite{gor} i.e., to demonstrate that an appropriate 
real-valued function on $X_M$, in our case 
$\log\varphi^{-1}:X_M\rightarrow\R$, is proper (i.e., the preimages of 
compact subsets are compact) with bounded 
gradient in modulus with respect to $g$ implying the completeness.        

Referring to (\ref{ujsokasag}) the open space $X_M$ arises by deleting 
one-one projective line from each $\C P^2$ factor, respectively, of 
the closed space $\overline{X}_M$. First we observe that 
$\varphi^{-1}:X_M\rightarrow\R_+$ is uniformly divergent along these 
projective lines. Assume that $\varphi^{-1}$ 
extends over $\overline{X}_M\supset X_M$ in a uniform continuous manner i.e. 
$\overline{\varphi}^{-1}\in C^0(\overline{X}_M)$ exists. A general 
principle based on the twistor construction is that the continuous 
extendibility of $\varphi^{-1}$ over $U\subseteqq \overline{X}_M$ implies the 
extendibility of $\pi$ i.e. the holomorphic map (\ref{celfuggveny}) over 
$\overline{p}^{-1}(U)\subseteqq\overline{Z}$ too in a manner which is 
compatible with the real structure on $\overline{p}^{-1}(U)$ i.e. this 
extension is a trivialization of the real bundle $S(\wedge^-U)$ 
(see the summary of twistor theory in Section \ref{two}). Therefore by our 
assumption the holomorphic map (\ref{celfuggveny}) extends over $\overline{Z}$ 
as well. However, since $\C P^2$ is not spin 
$\overline{X}_M=M\#\C P^2\#\dots\#\C P^2$ cannot be spin, too; consequently 
$S(\wedge^-\overline{X}_M)$ underlying the compact twistor space 
$\overline{Z}$ of $(\overline{X}_M,\overline{\gamma})$ is 
topologically not trivial hence its globally trivializing map 
(\ref{celfuggveny}) cannot extend from $X_M$ to $\overline{X}_M$, 
a contradiction. Assume that $\varphi^{-1}$ extends over 
at least one point of $\overline{X}_M\setminus X_M$ continuously. It yet  
follows that we run into a same type of contradiction. Assume now 
that $\varphi^{-1}$ extends over at least one point of 
$\overline{X}_M\setminus X_M$ in a discontinuous-but-bounded manner. Then we 
proceed as follows. The conformal scaling function satisfies with respect to 
$\gamma$ the following equations on $X_M$:
\begin{equation}
\left\{\begin{array}{ll}
\Delta_\gamma\varphi^{-1}+\frac{1}{6}\varphi^{-1}{\rm Scal}_\gamma&=0\:\:\:
\mbox{(vanishing of the scalar curvature of $g$ on $X_M$)};\\
&\\
\nabla^2_\gamma\varphi -\frac{1}{4}\left(\Delta_\gamma\varphi\right)
\gamma+\frac{1}{2}\varphi\:{\rm Ric}^0_\gamma&= 0\:\:\:
\mbox{(vanishing of the traceless Ricci tensor of $g$ on $X_M$)}.
\end{array}\right.
\label{skalazas}
\end{equation}
The Ricci tensor ${\rm Ric}_\gamma$ of $\gamma$ extends smoothly over 
$\overline{X}_M$ because it is just the restriction of the Ricci tensor of 
the self-dual metric $\overline{\gamma}$ on $\overline{X}_M$. Therefore 
both its scalar curvature ${\rm Scal}_\gamma$ and traceless Ricci part ${\rm 
Ric}^0_\gamma$ extend. Thus from the first equation of 
(\ref{skalazas}) we can see that 
$\varphi\Delta_\gamma\varphi^{-1}$ extends smoothly over $\overline{X}_M$. 
Likewise, adding 
the tracial part to the second equation of (\ref{skalazas}) we get 
$\varphi^{-1}\nabla^2_\gamma\varphi =-\frac{1}{2}{\rm Ric}_\gamma$ hence we 
conclude that the symmetric tensor field 
$\varphi^{-1}\nabla^2_\gamma\varphi$ extends 
smoothly over $\overline{X}_M$ so its trace 
$\varphi^{-1}\Delta_\gamma\varphi$ as well. 
Expanding $\Delta_{\overline{\gamma}}(\varphi\varphi^{-1})=0$ gives 
$(\Delta_{\overline{\gamma}}\varphi) \varphi^{-1} +2\:\overline{\gamma}
(\dd\varphi\:,\dd\varphi^{-1})+\varphi\Delta_{\overline{\gamma}}
\varphi^{-1}=0$ hence we obtain the pointwise equality 
\begin{equation}
\varphi^2\left\vert\dd\varphi^{-1}\right\vert^2_{\overline{\gamma}}=
\frac{1}{2}\left(\varphi\Delta_{\overline{\gamma}}\varphi^{-1}+
\varphi^{-1}\Delta_{\overline{\gamma}}\varphi\right)
\label{azonossag}
\end{equation} 
which demonstrates that $\varphi\left\vert\dd\varphi^{-1}\right\vert_{\overline
{\gamma}}$ extends smoothly over $\overline{X}_M$, too. If $\varphi^{-1}$ was 
extendible as a discontinuous bounded function over a point of 
$\overline{X}_M\setminus X_M$ then its gradient $\dd\varphi^{-1}$ was 
divergent in that point; hence from the extendibility of 
$\varphi\vert\dd\varphi^{-1}\vert_{\overline{\gamma}}$ we obtain that 
$\varphi$ vanishes hence $\varphi^{-1}$ is unbounded in that point, a 
contradiction again. We conclude that 
$\varphi^{-1}:X_M\rightarrow\R_+$ is {\it uniformly divergent along the 
whole complementum} $\overline{X}_M\setminus X_M$ yielding, on the one 
hand, that the function $\log\varphi^{-1}:X_M\rightarrow\R$ is proper. 
 
As a byproduct the inverse of $\varphi^{-1}$ is bounded on $X_M$ 
i.e., $\vert\varphi\vert\leqq c_1$ with a finite constant. 
We already know that $\vert\varphi\Delta_{\gamma}
\varphi^{-1}\vert\leqq c_2$ and $\vert\varphi^{-1}
\Delta_{\gamma}\varphi\vert\leqq c_3$ with other 
finite constants as well. 
Now writing $\varphi\vert\dd\varphi^{-1}\vert_{\gamma}=
\vert\dd (\log\varphi^{-1})\vert_{\gamma}$ and carefully noticing 
that $\vert\xi\vert_g=\varphi\vert\xi\vert_\gamma$ on
$1$-forms we can use (\ref{azonossag}) and the estimates above to come up with 
\[\vert\dd (\log \varphi^{-1})\vert_g
\leqq c_1\vert\dd (\log \varphi^{-1})\vert_{\gamma}
\leqq c_1\left(\left\vert\varphi\Delta_\gamma\varphi^{-1}\vert 
+\vert\varphi^{-1}\Delta_\gamma\varphi\right\vert\right)^{\frac{1}{2}}\leqq 
c_1(c_2+c_3)^{\frac{1}{2}}<+\infty\]
and conclude, on the other hand, that $\log\varphi^{-1}: 
X_M\rightarrow\R$ has bounded gradient in modulus with respect to $g$. 
Therefore, in light of Gordon's theorem \cite{gor}, the Ricci-flat space 
$(X_M,g)$ is complete. $\Diamond$
\vspace{0.1in}

\noindent We want to finish the construction by ending up with an 
open space with a single end, hence we want to remove the 
extra ``non-distinguished'' $\R^4$'s from $X_M$ in its decomposition 
(\ref{ujsokasag}) without destroying completeness and Ricci flatness. 
 
\begin{lemma} Consider the space $(X_M,g)$ as in Lemma 
\ref{teljeslemma}. Then the orientation and the complete Ricci-flat 
metric $g$ on $X_M$ descend to the punctured space $M^\times\subset M$ with 
its inherited smooth structure, rendering it a connected and simply connected, 
open, oriented, complete, Ricci-flat Riemannian spin $4$-manifold 
$(M^\times, g)$. 
\label{levagolemma}
\end{lemma}

\noindent {\it Proof}. It is clear from (\ref{ujsokasag}) 
that $M^\times$ arises from $X_M$ by filling in the ``centers'' of the 
finitely many non-distinguished $\R^4$ summands  with one-one point, 
respectively (see Figure 3). 
\vspace{0.4in}

\centerline{
\begin{tikzpicture}[scale=0.7]
\node at (10.2,0) {$M^\times$};
\draw [thick] plot [tension=0.8,smooth] coordinates 
{(11,0) (11.1,0.5) (11.5,1) (12.25,1.4) (13.5,1.5) (14.75,1.2) (16,0.2)};
\draw [thick] plot [tension=0.8,smooth] coordinates 
{(11,0) (11.1,-0.5) (11.5,-1) 
(12.25,-1.4) (13.5,-1.5) (14.75,-1.2) (16,-0.2)};
\draw [thick] plot [tension=0.8,smooth] coordinates {(16,0.2) (17,0.2) 
(18,1)};
\draw [thick] plot [tension=0.8,smooth] coordinates {(16,-0.2) (17,-0.2) 
(18,-1)};
\draw [gray] plot [tension=0.8,smooth] coordinates {(18,1) (18.3,0) 
(18,-1)};
\draw [gray] plot [tension=0.8,smooth] coordinates {(18,1) (17.7,0) 
(18,-1)};
\draw [thick] plot [tension=0.8,smooth] coordinates {(11.5,0.4) (12.5,0.2)
(13.5,0.4)};
\draw [thick] plot [tension=0.8,smooth] coordinates {(11.8,0.3) (12.5,0.5)
(13.2,0.3)};
\node at (-2.5,0) {$X_M$};
\draw [thick] plot [tension=0.8,smooth] coordinates { (0.5,0.2) (1,-0.2)
(1.5,0.2)};
\draw [thick] plot [tension=0.8,smooth] coordinates { (0.7,0) (1,0.1) (1.3,0)};
\draw [thick] plot [tension=0.8,smooth] coordinates
{(0.5,1) (1,1.3) (2,1) (3,0.7) (4,0.9) (4.5,1.3)};  
\draw [thick] plot [tension=0.8,smooth] coordinates
{(5.1,1.3) (5.4,0.6) (6,0.3)};
\draw plot [tension=0.8,smooth] coordinates
{(4.5,1.3) (4.8,1.4) (5.1,1.3)};
\draw plot [tension=0.8,smooth] coordinates
{(4.5,1.3) (4.8,1.2) (5.1,1.3)};
\draw [thick] plot [tension=0.8,smooth] coordinates 
{(0.5,-1) (1,-1.3) (2,-1) (3,-0.7) (4,-0.9) (4.5,-1.3)};                
\draw [thick] plot [tension=0.8,smooth] coordinates
{(5.1,-1.3) (5.4,-0.6) (6,-0.3)};
\draw plot [tension=0.8,smooth] coordinates
{(4.5,-1.3) (4.8,-1.4) (5.1,-1.3)};
\draw plot [tension=0.8,smooth] coordinates
{(4.5,-1.3) (4.8,-1.2) (5.1,-1.3)};
\draw [thick] plot [tension=0.8,smooth] coordinates {(0.5,1) (0,0) (0.5,-1)};
\draw [thick] plot [tension=0.8,smooth] coordinates {(6,0.3) (7,0.2) (8,1)};
\draw [thick] plot [tension=0.8,smooth] coordinates {(6,-0.3) (7,-0.2) (8,-1)};
\draw [gray] plot [tension=0.8,smooth] coordinates {(8,1) (7.8,0) (8,-1)};
\draw [gray] plot [tension=0.8,smooth] coordinates {(8,1) (8.2,0) (8,-1)};
%\label{abra3}
\end{tikzpicture}
}
\vspace{0.1in}

\centerline{Figure 3. Construction of $M^\times$ out of $X_M$ by 
filling in the extra $\R^4$'s.}
\vspace{0.3in}

\noindent Given this set-up, our strategy to prove the lemma is as follows: 
First apply Uhlenbeck's singularity removal theorem at each $\R^4$ summand to 
get rid of the corresponding singularity of the Levi--Civita connection of 
$g$---which is certainly an obstacle against the extension of the 
metric over the ``center'' of this $\R^4$ summand in the intermediate 
manifold---and in this way extend the connection to $M^\times$. Finally 
around each former singular point use a geodesic normal coordinate system 
adapted to this extended smooth connection on $M^\times$ to conclude that 
the metric $g$ on $X_M$ 
smoothly extends over the singularities, too. If this procedure works then 
the result is a smooth complete Ricci-flat metric on $M^\times$. However, 
as we shall see shortly, the non-existence of a spin structure on the 
original compact $M$ plays the role of an (and the only one) obstruction 
against the feasibility of this procedure. 

So let us take a fixed $\R^4$ summand in 
$X_M=M^\times\#\R^4\#\dots\#\R^4$. Since $X_M$ locally looks like a 
punctured $M$ around this summand i.e., a point $p\in M$ removed, we can 
{\it diffeomorphically} model $M^\times\#\R^4\#\dots\#\R^4$ around this 
$\R^4$ summand by an open punctured ball in some local modeling $\R^4$. 
More precisely let $p\in M$ be a point, $p\in U\subset M$ a neighbourhood 
containing the point and consider a local coordinate system 
$(U,y_1,\dots,y_4)$ centered at 
$p$ i.e., satisfying $y_1(p)=0,\dots, y_4(p)=0$. Identifying this local 
coordinate system with $(x_1,\dots,x_4)$ about the origin of the 
modeling $\R^4$ implies that $p$ is mapped to $0\in\R^4$ having coordinates 
$(x_1,\dots,x_4)=(0,\dots,0)$ and our model for the vicinity of 
the given $\R^4$ summand in $X_M$ then looks like 
\begin{equation}
\left( B^{4\times}_r(0)\:,\: x_1,\dots, x_4\right) 
\label{golyo}
\end{equation}
i.e., a coordinatized open punctured ball $B^{4\times}_r(0):=
B^4_r(0)\setminus\{0\}$ about $0\in\R^4$ of (Euclidean) radius $r>0$. (In 
this picture the ``infinity'' of the $\R^4$ summand corresponds to the 
the center of the ball.) Consider the restricted tangent bundle 
$TB^{4\times}_r(0):=TX_M\vert_{B^{4\times}_r(0)}$; using the 
restrictions of the orientation on $X_M$ and the metric $g$, we can render 
it a real four-rank ${\rm SO}(4)$ vector bundle over the punctured ball 
$B^{4\times}_r(0)$. We claim that $TB^{4\times}_r(0)$ in fact can be 
reduced to a complex two-rank ${\rm SU}(2)\subset{\rm SO}(4)$ vector bundle 
over the annulus. We can see this by exploiting the so far unmentioned 
feature of our construction namely that as a ``byproduct'' 
the space $(X_M,g)$ of Lemma \ref{teljeslemma} carries a compatible 
hyper-K\"ahler structure, too. Since the original compact space 
$(\overline{X}_M,\overline{\gamma})$ of Lemma \ref{kezdolemma} was 
oriented and self-dual with both properties being conformally invariant, 
$(X_M,g)$ is in fact a connected, simply connected, oriented, complete 
self-dual and Ricci-flat space or in other words: A hyper-K\"ahler 
$4$-manifold \cite[Chapter 13]{bes}. This implies among other things that 
the holonomy group of the Levi--Civita connection of $g$ hence the structure 
group of $TX_M$ reduces to ${\rm SU}(2)\subset {\rm SO}(4)$. Consider the 
Levi--Civita connection of $(X_M,g)$. We can therefore suppose that its 
restriction to $TB^{4\times}_r(0)\subset TX_M$ is 
an ${\rm SU}(2)$ connection $\nabla^\times$ suffering from a singularity at 
the origin. We know moreover that being $\nabla^\times$ self-dual, it solves 
the ${\rm SU}(2)$ Yang--Mills equations. Moreover $\nabla^\times$ has 
finite energy over $B^{4\times}_r(0)$. This is because $g$ is Ricci-flat 
and self-dual so the curvature of $\nabla^\times$ coincides with the 
self-dual Weyl component $W^+_g$ of $g$ only; however being conformally 
invariant, $W^+_g=W^+_\gamma=W^+_{\overline{\gamma}}\vert_{X_M}$ that is, 
the curvature tensor of $\nabla^\times$ is just the restriction of the Weyl 
tensor of the original smooth metric $\overline{\gamma}$ on 
$\overline{X}_M$. Consequently it is a smooth bounded tensor field on 
$B^{4\times}_r(0)$ implying finite local energy. This 
also yields that, if $0<r$ is sufficiently small, $\nabla^\times$ 
admits an $L^2_1$ gauge along $B^{4\times}_r(0)$ as well. Therefore, 
by Uhlenbeck's singularity removal theorem (see Theorem \ref{uhlenbecktetel}) 
there exists an $L^2_2$ gauge transformation on $TB^{4\times}_r(0)$ 
such that the gauge transformed connection extends across the singularity to 
a smooth ${\rm SU}(2)$ connection $\nabla$ on the trivial bundle $TB^4_r(0)$. 
Consequently, switching to the global picture, the singularity of the 
Levi--Civita connection around the fixed $\R^4$ summand of $X_M$ 
can be removed hence the corresponding $\R^4$ summand 
can be deleted from (\ref{ujsokasag}) according to our original plan. 
Repeating this procedure around all the finitely many $\R^4$ summands 
of $X_M$ we finally come up with a smooth ${\rm SU}(2)$ connection over 
$M^\times$. 

However there is an important topological subtlety here. For notational 
simplicity suppose that $X_M=M^\times\#\R^4$ i.e., possesses one 
non-distinguished $\R^4$ summand only. Then the singularity removal procedure 
carried out above convinces us that the original singular Levi--Civita 
connection defined on the {\it tangent} bundle $T(M^\times\#\R^4)$, regarded 
as an ${\rm SU}(2)$ bundle, indeed extends to a non-singular 
${\rm SU}(2)$ connection on {\it some} ${\rm SU}(2)$ bundle $E^\times$ over 
$M^\times$ i.e., it indeed smoothly exists somewhere which is however not 
necessarily the tangent bundle of $M^\times$. For instance, as we 
emphasized in the discussion after Theorem \ref{uhlenbecktetel}, the 
singularity-removing-gauge-transformation is not continuous in general 
hence the original global vector bundle carrying the singular connection may 
change topology during the singularity removal procedure. However, we know 
the following two things. On the one hand complex two-rank ${\rm SU}(2)$ 
vector bundles over $M^\times$, like the $E^\times$ above carrying the 
non-singular connection, are classified by various characteristic classes 
taking values in the groups $H^i(M^\times;\pi_{i-1}({\rm SU}(2))$ 
with $i=1,\dots,4$. Knowing the first three homotopy groups of ${\rm SU}(2)$ 
and taking into account the non-compactness of $M^\times$ these 
cohomology groups are all trivial consequently we know that $E^\times$ 
is necessarily isomorphic to the trivial bundle over $M^\times$. On the 
other hand, real rank-four ${\rm SO}(4)$ vector bundles over $M^\times$, like 
the tangent bundle $TM^\times$ carrying an orientation and a Riemann 
metric, are classified by characteristic classes taking 
values in $H^i(M^\times;\pi_{i-1}({\rm SO}(4))$. 
Again recalling the first three homotopy groups of the non-simply 
connected group ${\rm SO}(4)$ and still keeping in mind that $M^\times$ is 
non-compact, the only potentially non-trivial group here is 
$H^2(M^\times;\Z_2)$ demonstrating that vector bundles of this type over 
$M^\times$ are classified by a single element and this is nothing but their 
second Stiefel--Whitney class. Consequently if $M$ 
is spin or equivalently $w_2(TM)=0\in H^2(M;\Z_2)$ then by the injection 
$M^\times\subset M$ we find $w_2(TM^\times )=0
\in H^2(M^\times;\Z_2)$ as well showing that 
{\it $TM^\times$ is isomorphic to the trivial bundle}, too. Therefore we 
conclude that whenever $M$ is spin, we can identify the vector 
bundle $E^\times$ carrying the non-singular ${\rm SU}(2)$ connection over 
$M^\times$ with its tangent bundle $TM^\times$. 

Having understood this, we can finish the proof by extending the metric 
itself through the singularities. Fortunately this is simple. Consider the 
restricted connection $\nabla$ about one singular point $p$. This is now 
an overall (i.e., including the singular point) smooth connection. 
Therefore there exists a $\delta (p)>0$ such that we can 
suppose without loss of generality that the coordinate system 
(\ref{golyo}) we take about this singular point with 
$0<r<\delta (p)$ is a geodesic normal coordinate system 
with respect to $\nabla$. This implies that the Christoffel symbols 
$\nabla_{\partial_i}\partial_j=\sum\limits_k\Gamma^k_{ij}\partial_k$ all 
vanish in the center i.e., $\Gamma_{ij}^k(0,\dots,0)=0$ for all 
$i,j,k=1,\dots 4$. Then the well-known compatibility equations
\[\Gamma_{ij}^k=\frac{1}{2}
\sum\limits_{l=1}^4\left(\partial_ig_{lj}+\partial_jg_{li}-
\partial_lg_{ij}\right)g^{lk}\]
imply in a well-known way that in this gauge $g$ 
extends over the origin, too, such that $g_{ij}(0,\dots,0)=\delta_{ij}$ 
and $\partial_kg_{ij}(0,\dots,0)=0$ for all $i,j,k=1,\dots,4$. The further 
differentiablity i.e., the smoothness of $g$ at the origin follows from 
the smoothness of the Christoffel symbols there. That $g$ is 
Ricci-flat is a trivial consequence of the same property of the 
original metric. $\Diamond$

%%%%%%%%%%%%%%%%%%%%%%%%%%%%%%%%%%%

\section{Construction in the exotic setting}
\label{four}

%%%%%%%%%%%%%%%%%%%%%%%%%%%%%%%%%%%%%%%%%%%%

In this section we shall sink into the bottomless sea of four
dimensionality, called Exotica, and repeat the procedure performed in 
Section \ref{three}. That is, we shall construct solutions of 
the vacuum Einstein equation on the smooth $4$-manifold $X^\times$ which 
is only {\it homeomorphic} but not {\it diffeomorphic} to the punctured 
manifold $M^\times$ appeared in Section \ref{three}. This construction 
basically goes along the lines of that presented in Section \ref{three} 
with minor technical differences. Consequently, those steps which 
require new tools will be worked out in detail while those which are 
basically the same as the corresponding ones in Section \ref{three} will 
be sketched only.

To begin with, we compose Theorems \ref{taubestetel}, 
\ref{egzotikusnagycsalad} and \ref{gompftetel} together as follows.

\begin{lemma}

Out of any connected, closed (i.e., compact without boundary) oriented
smooth $4$-manifold $M$ one can construct a connected, open (i.e.,
non-compact without boundary) oriented smooth Riemannian $4$-manifold
$(X_M,\gamma)$ which is self-dual but incomplete in general.
\label{egzotikuskezdolemma}
\end{lemma}

\noindent {\it Proof.} Pick any connected, oriented, closed, smooth
$4$-manifold $M$. Referring to Theorem \ref{taubestetel} let 
$k:=\max (1,k_M)\in\N$ be a positive integer, put
\[\overline{X}_M:=M\#\underbrace{\C P^2\#\dots\#\C P^2}_{k}\]
and let $\overline{\gamma}$ be a self-dual metric on it. Then
$(\overline{X}_M,\overline{\gamma})$ is a compact self-dual manifold. Pick 
one $\C P^2$ factor within $\overline{X}_M$ and denote by 
$S^2:=\C P^2\setminus R^4$ the complement of the largest exotic 
$\R^4$-space $R^4\subset\C P^2$, considered as an only 
``continously embedded projective line'' in that factor, as in the 
discussion after Theorem \ref{egzotikusnagycsalad} (we can suppose that the 
closed subspace $S^2\subset\C P^2$ avoids the attaching point of $\C P^2$ to 
$M$). Let $K\subset R^4$ be the connected compact subset as in part (ii) of 
Theorem \ref{egzotikusnagycsalad} and put
\begin{equation}
X_M:=M\#\underbrace{(\C P^2\setminus\C P^1)\#\dots\#
(\C P^2\setminus\C P^1)}_{k-1}\#_K(\C P^2\setminus S^2)\cong 
M\#\underbrace{\R^4\#\dots\#\R^4}_{k-1}\#_KR^4\cong 
X^\times\#\underbrace{\R^4\#\dots\#\R^4}_{k-1}
\label{egzotikusujsokasag}
\end{equation}
where the operation $\#_K$ means that the attaching point $y_0\in R^4$
used to glue $R^4$ with $M\#\R^4\#\dots\#\R^4$ satisfies
$y_0\in K\subset R^4$ and $X^\times:=M\#_KR^4$ is a smooth manifold 
homeomorphic but not diffeomorphic to the puncturation $M^\times$ of the 
original manifold (see Theorem \ref{gompftetel}). The result is a 
connected, open $4$-manifold $X_M$ (see Figure 4). 
\vspace{0.3in}

\centerline{
\begin{tikzpicture}[scale=0.7]
\node at (-0.8,0) {$M$};
\draw [thick] plot [tension=0.8,smooth cycle] coordinates {
(0,0) (2,2) (5,0) (2,-2)};
\draw [thick] plot [tension=0.8,smooth] coordinates {(0.5,0.4) (1.5,0.2)
(2.5,0.4)};
\draw [thick] plot [tension=0.8,smooth] coordinates {(0.8,0.3) (1.5,0.5)
(2.2,0.3)};
\node at (7.5,0) {$X_M$};
\draw [thick] plot [tension=0.8,smooth] coordinates { (9.5,0.2) 
(10.2,-0.2) (10.9,0.2)};
\draw [thick] plot [tension=0.8,smooth] coordinates { (9.7,0) 
(10.2,0.3) (10.7,0)};
\draw [thick] plot [tension=0.8,smooth] coordinates
{(9.5,1) (10,1.3) (11,1) (12,0.7) (13,0.9) (13.5,1.3)};
\draw [thick] plot [tension=0.8,smooth] coordinates
{(14.1,1.3) (14.4,0.6) (15,0.3)};
\draw plot [tension=0.8,smooth] coordinates
{(13.5,1.3) (13.8,1.4) (14.1,1.3)};
\draw plot [tension=0.8,smooth] coordinates
{(13.5,1.3) (13.8,1.2) (14.1,1.3)};
\draw [thick] plot [tension=0.8,smooth] coordinates
{(9.5,-1) (10,-1.3) (11,-1) (12,-0.7) (13,-0.9) (13.5,-1.3)};
\draw [thick] plot [tension=0.8,smooth] coordinates
{(14.1,-1.3) (14.4,-0.6) (15,-0.3)};
\draw plot [tension=0.8,smooth] coordinates
{(13.5,-1.3) (13.8,-1.4) (14.1,-1.3)};
\draw plot [tension=0.8,smooth] coordinates
{(13.5,-1.3) (13.8,-1.2) (14.1,-1.3)};
\draw [thick] plot [tension=0.8,smooth] coordinates {(9.5,1) (9,0)
(9.5,-1)};
\draw [thick] plot [tension=0.8,smooth] coordinates {(15,0.3) (16,0.2) 
(17,1)};
\draw [thick] plot [tension=0.8,smooth] coordinates {(15,-0.3) (16,-0.2)
(17,-1)};
\draw[gray] plot coordinates {(17,1) (17.5,0.9) (17,0.7) (17.5,0.6)
(17,0.5) (17.5,0.4) (17,0.3) (17.5,0.2) (17,0.1) (17.5,0) (17,-0.1)
(17.5,-0.2) (17,-0.3) (17.5,-0.4) (17,-0.5) (17.5,-0.6) (17,-0.7)
(17.5,-0.9) (17,-1)};
%\label{abra4}
\end{tikzpicture}
}
\vspace{0.1in}

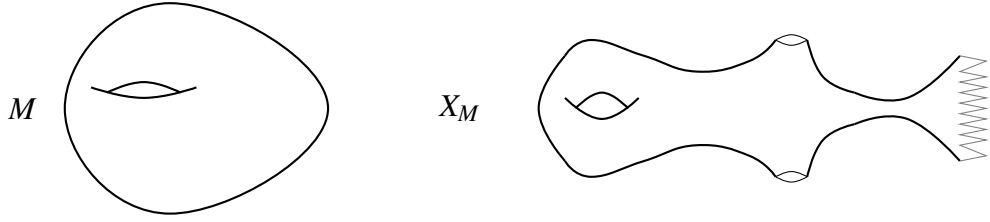
\centerline{Figure 4. Construction of $X_M$ out of $M$ in the exotic 
setting. The gray zig-zag represents a}  
\centerline{``creased end'' diffeomorphic to the complement of a connected 
compact subset $K$ in the exotic $R^4$.}
\vspace{0.3in}

\noindent From the proper smooth embedding $X_M\subsetneqq\overline{X}_M$ 
there exists a restricted self-dual Riemannian metric
$\gamma:=\overline{\gamma}\vert_{X_M}$ on $X_M$ which is however in
general non-complete. $\Diamond$
\vspace{0.1in}

\noindent In the case of our situation set up in Lemma 
\ref{egzotikuskezdolemma}
twistor theory works as follows. Consider the compact self-dual
space $(\overline{X}_M,\overline{\gamma})$ from Lemma 
\ref{egzotikuskezdolemma}, take its twistor fibration $\overline{p}: 
\overline{Z}\rightarrow\overline{X}_M$ and let $p:Z\rightarrow X_M$ 
be its restriction induced by the smooth embedding $X_M\subsetneqq 
\overline{X}_M$ i.e., $Z:=\overline{Z}\vert_{X_M}$ and
$p:=\overline{p}\vert_{X_M}$. Then $Z$ is a non-compact
complex $3$-manifold already obviously possessing all the required twistor  
data except the existence of a holomorphic mapping $\pi :Z\rightarrow\C P^1$.

\begin{lemma} Consider the connected, open, oriented, incomplete, self-dual 
space $(X_M,\gamma)$ as in Lemma 
\ref{egzotikuskezdolemma} with its twistor fibration $p:Z\rightarrow X_M$ 
constructed above. If $\pi_1(M)=1$ and $M$ is spin (or equivalently, having 
even intersection form) then there exists a holomorphic mapping 
$\pi :Z\rightarrow\C P^1$.

\label{egzotikuspilemma}
\end{lemma}

\noindent {\it Proof.} Let $x_0\in X_M$ be an arbitrary fixed point 
of $X_M$ in (\ref{egzotikusujsokasag}). Our aim is to construct a 
holomorphic map
\[\pi:\:\:Z\longrightarrow p^{-1}(x_0)\cong\C P^1\]
that we carry out exactly the same way as in the proof of Lemma 
\ref{pilemma} hence we do not repeat it here. $\Diamond$ 
\vspace{0.1in}

\noindent It also follows that $\pi :Z\rightarrow\C P^1$ i.e., the map 
constructed in Lemma \ref{egzotikuspilemma} is 
compatible with the real structure $\tau :Z\rightarrow Z$ already fixed by 
the self-dual structure in Theorem \ref{taubestetel} therefore twistor theory 
provides us with a Ricci-flat (and self-dual) Riemannian metric $g$ 
on $X_M$. We proceed further and demonstrate that, unlike 
$(X_M,\gamma)$, the space $(X_M,g)$ is complete.

\begin{lemma} The four dimensional connected and simply connected, open, 
oriented, Ricci-flat Riemannian spin manifold $(X_M,g)$ is complete. 
\label{egzotikusteljeslemma}
\end{lemma}

\noindent {\it Proof}. The metrics $\gamma$ and $g$ originate from the same 
twistor space again hence they are conformally equivalent consequently there 
exists a smooth function 
$\psi :X_M\rightarrow\R_+$ satisfying $g=\psi^{-2}\gamma$. 
Taking into account that the steps in the proof of Lemma \ref{teljeslemma} 
have been insensitive for the particular construction of the complementum 
$\overline{X}_M\setminus X_M$ we can simply repeat them here. Hence 
we find again, on the one hand, that $\psi^{-1}$ blows up 
uniformly along the whole $\overline{X}_M\setminus X_M$ this time consisting of 
the disjoint union of ``ordinary'' i.e. holomorphically embedded 
projective lines $\C P^1=\C P^2\setminus\R^4$ and the distinguished 
``continuously embedded projective line'' $S^2=\C P^2\setminus R^4$ 
in the distinguished factor in (\ref{egzotikusujsokasag}); consequently 
$\log\psi^{-1}$ is proper. Moreover, on the other hand, recalling the 
steps of Lemma \ref{teljeslemma} we see that $\log\psi^{-1}$ has 
bounded gradient in modulus with respect to $g$, too. Consequently $(X_M, g)$ 
is complete as in the proof of Lemma \ref{teljeslemma} hence the details are 
omitted. $\Diamond$

\begin{remark}\rm For clarity we remark that comparing the proofs of 
Lemmata \ref{teljeslemma} and \ref{egzotikusteljeslemma} one cannot 
conclude that the resulting complete spaces $(X_M,g)$ in the 
non-exotic and exotic situations are conformally equivalent. This is because 
(see the discussion after Theorem \ref{egzotikusnagycsalad}) 
the locations of $X_M$ in the two cases within their common 
closure $\overline{X}_M$ are different such that even the former cannot be 
mapped into the latter by any diffeomorphism of $\overline{X}_M$. Consequently 
taking the pointwise product of the scaling function $\varphi^{-2}$ in 
Lemma \ref{teljeslemma} with the inverse one $\psi^2$ from Lemma 
\ref{egzotikusteljeslemma} to obtain a conformal rescaling 
between the corresponding metrics makes no sense. 
\end{remark}

\noindent Finally we cut down the standard $\R^4$'s from $X_M$ to 
obtain $X^\times$ as in Lemma \ref{levagolemma}.

\begin{lemma} Consider the space $(X_M,g)$ as in Lemma
\ref{egzotikusteljeslemma}. Then the orientation and the complete 
Ricci-flat metric $g$ on $X_M$ descend to the punctured space $X^\times$ 
(which is homeomorphic but not diffeomorphic to the corresponding space 
$M^\times$ of Lemma \ref{levagolemma}) with
its inherited smooth structure, rendering it a connected and simply connected,
open, oriented, complete, Ricci-flat Riemannian spin $4$-manifold
$(X^\times, g)$.
\label{egzotikuslevagolemma}
\end{lemma}

\noindent {\it Proof.} Taking into account that filling in the standard 
$\R^4$'s in the decomposition (\ref{egzotikusujsokasag}) of $X_M$ is a 
completely local procedure (see Figure 5)
\vspace{0.3in}

\centerline{
\begin{tikzpicture}[scale=0.7]
\node at (10.2,0) {$X^\times$};
\draw [thick] plot [tension=0.8,smooth] coordinates
{(11,0) (11.1,0.5) (11.5,1) (12.25,1.4) (13.5,1.5) (14.75,1.2) (16,0.2)};
\draw [thick] plot [tension=0.8,smooth] coordinates
{(11,0) (11.1,-0.5) (11.5,-1)
(12.25,-1.4) (13.5,-1.5) (14.75,-1.2) (16,-0.2)};
\draw [thick] plot [tension=0.8,smooth] coordinates {(16,0.2) (17,0.2)
(18,1)};
\draw [thick] plot [tension=0.8,smooth] coordinates {(16,-0.2) (17,-0.2)
(18,-1)};
\draw[gray] plot coordinates {(18,1) (18.5,0.9) (18,0.7) (18.5,0.6)
(18,0.5) (18.5,0.4) (18,0.3) (18.5,0.2) (18,0.1) (18.5,0) (18,-0.1)
(18.5,-0.2) (18,-0.3) (18.5,-0.4) (18,-0.5) (18.5,-0.6) (18,-0.7)
(18.5,-0.9) (18,-1)};
\draw [thick] plot [tension=0.8,smooth] coordinates {(11.5,0.4) (12.5,0.2)
(13.5,0.4)};
\draw [thick] plot [tension=0.8,smooth] coordinates {(11.8,0.3) (12.5,0.5)
(13.2,0.3)};
\node at (-1,0) {$X_M$};
\draw [thick] plot [tension=0.8,smooth] coordinates { (0.5,0.2) (1,-0.2)
(1.5,0.2)};
\draw [thick] plot [tension=0.8,smooth] coordinates { (0.7,0) (1,0.1) 
(1.3,0)};
\draw [thick] plot [tension=0.8,smooth] coordinates
{(0.5,1) (1,1.3) (2,1) (3,0.7) (4,0.9) (4.5,1.3)};
\draw [thick] plot [tension=0.8,smooth] coordinates
{(5.1,1.3) (5.4,0.6) (6,0.3)};
\draw plot [tension=0.8,smooth] coordinates
{(4.5,1.3) (4.8,1.4) (5.1,1.3)};
\draw plot [tension=0.8,smooth] coordinates
{(4.5,1.3) (4.8,1.2) (5.1,1.3)};
\draw [thick] plot [tension=0.8,smooth] coordinates
{(0.5,-1) (1,-1.3) (2,-1) (3,-0.7) (4,-0.9) (4.5,-1.3)};
\draw [thick] plot [tension=0.8,smooth] coordinates
{(5.1,-1.3) (5.4,-0.6) (6,-0.3)};
\draw plot [tension=0.8,smooth] coordinates
{(4.5,-1.3) (4.8,-1.4) (5.1,-1.3)};
\draw plot [tension=0.8,smooth] coordinates
{(4.5,-1.3) (4.8,-1.2) (5.1,-1.3)};
\draw [thick] plot [tension=0.8,smooth] coordinates {(0.5,1) (0,0)
(0.5,-1)};
\draw [thick] plot [tension=0.8,smooth] coordinates {(6,0.3) (7,0.2) (8,1)};
\draw [thick] plot [tension=0.8,smooth] coordinates {(6,-0.3) (7,-0.2) 
(8,-1)};
\draw[gray] plot coordinates {(8,1) (8.5,0.9) (8,0.7) (8.5,0.6)
(8,0.5) (8.5,0.4) (8,0.3) (8.5,0.2) (8,0.1) (8.5,0) (8,-0.1)
(8.5,-0.2) (8,-0.3) (8.5,-0.4) (8,-0.5) (8.5,-0.6) (8,-0.7)
(8.5,-0.9) (8,-1)};
%\label{abra5}
\end{tikzpicture}
}
\vspace{0.1in}

\centerline{Figure 5. Construction of $X^\times$ out of $X_M$ by filling 
in the extra $\R^4$'s.}
\vspace{0.3in}

\noindent the proof of this lemma is verbatim the same as the proof of Lemma 
\ref{levagolemma} hence is omitted. $\Diamond$
\vspace{0.1in}

\noindent{\it Proof of Theorems \ref{fotetel1} or \ref{fotetel2}}. 
Collecting all the ruslts of Sections \ref{three} and 
\ref{four} together the desired statements are obtained. $\Diamond$

\begin{remark}\rm Before proceeding further let us note that Theorems 
\ref{fotetel1} or \ref{fotetel2} correspond to the case when the 
creased end of $X^\times$ is diffeomorphic to the largest member 
$R^4=R^4_{+\infty}$ in the radial family of Theorem 
\ref{egzotikusnagycsalad} i.e. $X^\times=M\#R^4_{+\infty}$. It would be 
interesting to understand whether or not a similar construction works for 
the intermediate members of the family i.e. for $X^\times_t=M\#R^4_t$. 
\end{remark}

%%%%%%%%%%%%%%%%%%%%%%%%%%%%%%%%%%%%%%%%

\section{Lorentzian solutions}
\label{five}

%%%%%%%%%%%%%%%%%%%%%%%%%%%%%%%%%%%%%%

In the previous sections we have produced an immense class of Ricci flat 
Riemannian spaces $(X^\times ,g)$ which are non-compact but complete. In 
this section we convert all of them i.e. the spaces in Theorem \ref{fotetel1} 
or equivalently Theorem \ref{fotetel2} into Ricci-flat 
Lorentzian ones as formulated in Theorem \ref{lorentztetel} by 
(essentially verbatim) recalling \cite[Lemma 4.2]{ete2}. Conversion is 
in principle possible because all the underlying manifolds 
$X^\times$ are non-compact hence there is no topological obstruction 
against Lorentzian structure.
\vspace{0.1in}

\noindent {\it Proof of Theorem \ref{lorentztetel}}. By virtue of its global 
triviality (cf. Lemma \ref{levagolemma} or \ref{egzotikuslevagolemma}), 
$TX^\times$ admits a nowhere vanishing smooth section 
yielding a splitting $TX^\times=L\oplus L^\perp$ into a real line bundle 
$L\subset TX^\times$ spanned the section and its $g$-orthogonal complementum 
subbundle $L^\perp\subset TX^\times$. Take the 
complexification $T^\C X^\times:=TX^\times\otimes_\R\C$ of the real tangent 
bundle as well as the complex bilinear extension of the Riemannian Ricci-flat 
metric $g$ found on $TX^\times$ to a Ricci-flat metric $g^\C$ on 
$T^\C X^\times$. This means that if $v^\C$ is a complexified tangent vector 
then both $v^\C\mapsto g^\C(v^\C\:,\:\cdot\:):=g(v^\C\:,\:\cdot\:)$ and 
$v^\C\mapsto g^\C(\:\cdot\:,\:v^\C):=g(\:\cdot\:,\:v^\C)$ are declared to be 
$\C$-linear and of course ${\rm Ric}_{g^\C}={\rm Ric}_g=0$. There is an 
induced splitting 
\begin{equation}
T^\C X^\times=L\oplus L^\perp\oplus\sqrt{-1}\:L\oplus\sqrt{-1}\:L^\perp
\label{felbontas}
\end{equation}
over $\R$ of the complexification i.e., if $T^\C X^\times$ is 
considered as a real rank-$8$ bundle over $X^\times$. Define a metric on the 
real rank-$4$ sub-bundle $L^\perp\oplus\sqrt{-1}\:L\subset T^\C X^\times$ by 
taking the restriction $g^\C\vert_{L^\perp\oplus\sqrt{-1}\:L}$. It readily 
follows from the orthogonality and reality of the splitting that this is 
a non-degenerate real-valued $\R$-bilinear form of Lorentzian type on this 
real sub-bundle. To see this, we simply have to observe that taking real vector 
fields $v_1,v_2:X^\times\rightarrow L$ and $w_1,w_2:X^\times\rightarrow 
L^\perp$ we can exploit the $\C$-bilinearity of $g^\C$ to write 
\[g^\C\vert_{L^\perp\oplus\sqrt{-1}\:L}(\sqrt{-1}\:v_1,\sqrt{-1}\:v_1)=
g^\C(\sqrt{-1}\:v_1,\sqrt{-1}\:v_1)=-g^\C(v_1,v_1) =-g(v_1,v_1)\]
and
\[g^\C\vert_{L^\perp\oplus\sqrt{-1}\:L}(\sqrt{-1}\:v_1,w_1)=
g^\C(\sqrt{-1}\:v_1,w_1)=\sqrt{-1}\:g^\C(v_1,w_1)
=\sqrt{-1}g(v_1,w_1)=0\]
and finally
\[g^\C\vert_{L^\perp\oplus\sqrt{-1}\:L}(w_1,w_2)=
g^\C(w_1,w_2)=g(w_1,w_2)\:\:\:.\]
Consider the $\R$-linear bundle isomorphism $W_L:T^\C X^\times
\rightarrow T^\C X^\times$ of the complexified tangent bundle defined by, with
respect to the splitting (\ref{felbontas}), as
\[W_L(v_1,w_1,\sqrt{-1}\:v_2,\sqrt{-1}\:w_2):=
(v_2, w_1,\sqrt{-1}\:v_1,\sqrt{-1}\:w_2)\:\:.\] 
Obviously $W_L^2={\rm Id}_{T^\C X^\times}$ or more precisely $W_L$ is a 
{\it real reflection} with respect to $g^\C$ making the diagram 
\[\xymatrix{
\ar[d] T^\C X^\times\ar[r]^{W_L} & T^\C X^\times\ar[d] \\
          X^\times\ar[r]^{{\rm Id}_{X^\times}}  & X^\times
}\]
commutative. In particular it maps the real tangent bundle 
$TX^\times=L\oplus L^\perp\subset T^\C X^\times$ onto the real bundle
$L^\perp\oplus\sqrt{-1}\:L\subset T^\C X^\times$ and {\it vice versa}. 
Consequently with arbitrary two tangent vectors $v,w:X^\times\rightarrow 
TX^\times$ 
\[g_L(v,w):=g^\C(W_Lv\:,\:W_Lw)\]
satisfies $g_L(v,w)=g^\C\vert_{L^\perp\oplus\sqrt{-1}\:L}(W_Lv\:,\:W_Lw)$ 
i.e., obtain a non-degenerate real-valued $\R$-bilinear form of Lorentzian type 
hence a smooth Lorentzian metric $g_L$ on the original real tangent bundle 
$TX^\times$. 

Concerning the Ricci tensor of $g_L$, the Levi--Civita
connections $\nabla^L$ of $g_L$ and $\nabla^\C$ of $g^\C$ satisfy 
\begin{eqnarray}
g_L(\nabla^L_uv,w)+g_L(v,\nabla^L_uw)&=&\dd g_L(v,w)u\nonumber\\
&=&\dd g^\C(W_Lv,W_Lw)u\nonumber\\
&=&g^\C(\nabla^\C_u(W_Lv)\:,\:W_Lw)+g^\C(W_Lv\:,\:\nabla^\C_u(W_Lw))\nonumber\\
&=&g^\C(W^2_L\nabla^\C_u W_Lv\:,\:W_Lw)+g^\C(W_Lv\:,\:W^2_L\nabla^\C_uW_Lw)
\nonumber\\
&=&g_L((W_L\nabla^\C_u W_L)v\:,\:w)+g_L(v\:,\:(W_L\nabla^\C_u W_L)w)\nonumber
\end{eqnarray}
yielding $\nabla^L=W_L\nabla^\C W_L$ (this is an $\R$-linear operator). 
Consequently the curvature ${\rm Riem}_{g_L}$ of $g_L$ is 
\[{\rm Riem}_{g_L}(v,w)u=\left[\nabla^L_v\:,\nabla^L_w\right]u-
\nabla^L_{[v,w]}u=W_L({\rm Riem}_{g^\C}(v,w)W_Lu)\:\:.\]
Let $\{e_0,e_1,e_2,e_3\}$ be a real orthonormal frame for $g_L$ at 
$T_pX^\times$ satisfying $g_L(e_0,e_0)=-1$ and $+1$ for the rest; then
$W_Le_0=\sqrt{-1}\:e_0$ and $W_Le_j=e_j$ for $j=1,2,3$ together
with the definition of $g_L$ imply that
\[g_L({\rm Riem}_{g_L}(e_0,v)w,e_0)=g^\C\left(W_L({\rm
Riem}_{g_L}(e_0,v)w),W_Le_0\right)=g^\C\left({\rm
Riem}_{g^\C}(e_0,v)W_Lw\:,\sqrt{-1}e_0\right)\]
and likewise
\[g_L({\rm Riem}_{g_L}(e_j,v)w\:,\:e_j)=g^\C(W_L({\rm
Riem}_{g_L}(e_j,v)w),W_Le_j)=g^\C({\rm
Riem}_{g^\C}\left(e_j\:,\:v\right)W_Lw\:,e_j)\:\:.\]
Using an orthonormal frame $\{f_1,\dots,f_m\}$ for a metric $h$ of any 
signature, its Ricci tensor looks like ${\rm Ric}_h(v,w)=
\sum\limits_{k=1}^mh(f_k,f_k)h({\rm Riem}_h(f_k,v)w\:,\:f_k)$; hence
\begin{eqnarray}
{\rm Ric}_{g_L}(v,w)\!\!\!\!\!\!&=&\!\!\!\!\!\!g_L(e_0,e_0)g_L({\rm
Riem}_{g_L}(e_0,v)w,e_0)+
\sum\limits_{j=1}^3g_L(e_j,e_j)g_L({\rm Riem}_{g_L}(e_j,v)w,e_j)\nonumber\\
\!\!\!\!\!\!&=&\!\!\!\!\!\!g^\C(\sqrt{-1}e_0,\!\sqrt{-1}e_0)g^\C({\rm
Riem}_{g^\C}(e_0,v)W_Lw,\sqrt{-1}e_0)
\!+\!\!\sum\limits_{j=1}^3g^\C(e_j,e_j)g^\C({\rm
Riem}_{g^\C}(\!e_j,v)W_Lw,e_j\!)\nonumber\\
\!\!\!&=&\!\!\!(-\sqrt{-1}\:-1)g^\C(e_0,e_0)g^\C({\rm
Riem}_{g^\C}(e_0,v)W_Lw,e_0)+{\rm Ric}_{g^\C}(v\:,\:W_Lw)\nonumber\\
&=& (-1+\sqrt{-1})g_L({\rm Riem}_{g_L}(e_0,v)w,e_0)\nonumber
\end{eqnarray}
and we also used $\{e_0,e_1,e_2,e_3\}$ as a complex orthonormal 
basis for $g^\C$ on $T^\C_pX^\times$ to write 
\[\sum_{j=0}^3g^\C(e_j,e_j)g^\C({\rm Riem}_{g^\C}(e_j,v)W_Lw,e_j)=
{\rm Ric}_{g^\C}(v,W_Lw)=0\:\:.\] 
Being the left hand side in 
${\rm Ric}_{g_L}(v,w)=(-1+\sqrt{-1})g_L({\rm Riem}_{g_L}(e_0,v)w,e_0)$ real, 
the right hand side must be real as well for all $v,w\in T_pX^\times$ 
which is possible if and only if both sides vanish. This demonstrates that 
$g_L$ is indeed Ricci-flat. $\Diamond$

%%%%%%%%%%%%%%%%%%%%%%%%%%%%%%%%%

\section{Physical interpretation}
\label{six}

%%%%%%%%%%%%%%%%%%%%%%%%%%%%%%%%%%%%

In this closing section we discuss the physical interpretation of the 
Lorentzian Ricci-flat geometries found in Theorem \ref{lorentztetel}. We 
believe that an interpretation is necessary because there are many known 
physically irrelevant solutions of the vacuum Einstein equation and our 
solutions as presented in Theorem \ref{lorentztetel} are admittedly very 
implicit and transcendental hence their physical significance, if any, is 
unclear yet. The offered interpretation fits well into the context of the 
celebrated {\it strong cosmic censorship conjecture} in its usual broad 
formulation (${\bf SCCC}$ for short) which is a hot topic recently (far 
from being complete cf. \cite{car-cos-des-hin-jan1, 
car-cos-des-hin-jan2,daf-luk,dia-epe-rea-san,dia-rea-san,ete1,ete2, 
ge-jia-wan-zha-zho,hod,lun-zie-car-cos-nat,mo-tia-wan-zha-zho}; for 
historical accounts see \cite{cha,ise,pen2}) and the so far hypothetical {\it 
topology changing} phenomena (again far from being 
complete, cf. e.g. \cite{dow,gib,hor}).

The current situation of the ${\bf SCCC}$ can perhaps be best 
summarized as a {\it puzzling dichotomy}: although there are some signs or 
hints for its (in)validity in {\it physically relevant} situations 
(like various black holes in asymptotically flat or de Sitter space-times 
filled with vacuum or various matter fields, etc. \cite{car-cos-des-hin-jan1, 
car-cos-des-hin-jan2,daf-luk,dia-epe-rea-san,dia-rea-san,ge-jia-wan-zha-zho,hod,
lun-zie-car-cos-nat,mo-tia-wan-zha-zho}), these 
are still not sharp enough to decide the status of the ${\bf SCCC}$ in these 
important cases. On the other hand there exists an superabundance of 
``exotic'' smooth solutions in which the ${\bf SCCC}$ 
clearly fails \cite{ete1,ete2} (namely the ones exhibited in 
Theorem \ref{lorentztetel}) however the physical meaning of these quite 
{\it purely mathematical} solutions is not clear yet. The reason for this 
latter issue is that, although being smooth solutions of the vacuum 
Einstein equation hence apparently relevant, the ${\bf SCCC}$ violating 
properties of these ``exotic'' solutions rest neither on some physical 
phenomenon nor on standard analytico-geometric properties of Lorentzian 
metrics; but rather based on subtle novel differentio-topological features 
(often called {\it exotica}) of four dimensional manifolds which have 
gradually been recognized in the underlying mathematical model of physical 
space-times from the 1980's onwards (cf. \cite[Introduction]{ete2}). Despite 
that no {\it a priori} principle has been introduced so 
far to exclude these curious and apparently fundamental mathematical 
discoveries from the game, they have not found their right places in 
theoretical physics yet \cite{ass-bra}.

The aim of this section is an effort to fill in this gap by offering a 
plausible and simple physical interpretation of the new ${\bf SCCC}$ 
violating solutions \cite{ete1,ete2} (i.e. the spaces exhibited in Theorem 
\ref{lorentztetel}). As an interesting observation it 
will turn out that, meanwhile the aforementioned classical situations in 
which ${\bf SCCC}$ breakdown has been examined belong to the well known 
static or stationary regime of general relativity, the new ${\bf SCCC}$ 
violating solutions are related with the yet unexplored deep 
dynamical regime of general relativity describing spatial topology changes as 
will be explained shortly. We also find that this dynamics appears as a 
{\it cosmologial redshift} for late time internal observers within these 
space-times. Therefore, quite unsurprisingly, one is tempted to say that as 
one moves from the static towards the dynamical regime, ${\bf SCCC}$ violating 
phenomena become more and more relevant in general relativity. 

Take any connected, simply connected, closed spin $4$-manifold $M$ and 
form the connected sum $X^\times:=M\#R^4$ as before (see Figure 5). 
It is easy to see (cf. the summary of the exotic stuff 
in Section \ref{two}) that $X^\times$ is homeomorphic to the 
punctured space $M^\times =M\setminus\{{\rm point}\}$ however cannot be 
diffeomorphic to it (with its usual inherited smooth 
structure from the smooth embedding $M^\times\subset M$) since $M^\times$ is 
diffeomorphic to $M\#\R^4$ meanwhile $X^\times$ by construction is 
diffeomorphic to $M\# R^4$ hence the ends of the two open spaces, although 
homeomorphic, are not diffeomorphic. Actually, from a general viewpoint, 
the appearance of non-compact $4$-manifolds carrying smooth structures like 
$X^\times$ i.e. having a ``creased end'' is much more typical. 
Theorem \ref{lorentztetel} then says that 
$X^\times$ always carries a Ricci-flat Lorentzian metric $g_L$. Having 
$X^\times$ a creased end implies that it surely cannot be written as a smooth 
product $\Sigma\times\R$ where $\Sigma$ is a $3$-manifold and $\R$ is the real 
line (with their unique smooth structures); however the existence of such a 
smooth splitting is a necessary condition of global hyperbolicity 
\cite{ber-san}. Consequently we arrive at a sort of heavy breakdown of 
the ${\bf SCCC}$ (in its usual broad formulation, cf. e.g. 
\cite{ete1, ete2}), namely
\vspace{0.1in}

\noindent $\overline{{\bf SCCC}}$. {\it The smooth Ricci-flat Lorentzian 
$4$-manifold $(X^\times ,g_L)$ in Theorem \ref{lorentztetel} is not 
globally hyperbolic and no (sufficiently large in an appropriate topological 
sense) perturbation of it can be globally hyperbolic.}
\vspace{0.1in}

\noindent Furthermore, Theorem \ref{fotetel1} says that $X^\times$ 
carries a complete Ricci-flat Riemannian metric $g$, too. As a by-product 
of the construction we have seen that fixing an appropriate 
orientation on $X^\times$ the metric $g$ is self-dual, too. However a 
simply-connected, complete Riemannian $4$-manifold which is both 
Ricci-flat and self-dual is in fact, as formulated in Theorem \ref{fotetel2}, 
hyper-K\"ahler (cf. e.g. \cite[Chapter 13]{bes}). Physically speaking the 
Riemannian $4$-manifolds $(X^\times, g)$ exhibited in Theorem \ref{fotetel1} 
or equiavlently, in Theorem \ref{fotetel2} are therefore examples of 
{\it gravitational instantons}. Consequently, even if these Riemannian 
(or Euclidean) vacuum spaces might not play any role in classical general 
relativity, they are not negligable in any quantum theory perhaps 
lurking behind classical general relativity.

After these introductory or general remarks let us move towards 
a suggested physical interpretation of the ${\bf SCCC}$ breaking 
but otherwise regular geometry $(X^\times, g_L)$. The conversion procedure 
in Theorem \ref{lorentztetel} rests on a nowhere-vanishing vector field 
\begin{equation}
v\in C^\infty (X^\times; TX^\times\setminus\{0\})
\label{vektormezo}
\end{equation}
along $X^\times$ whose choice was otherwise arbitrary. 
Therefore, taking into account the global triviality of the tangent bundle 
$TX^\times$ (cf. Lemmata \ref{levagolemma} and \ref{egzotikuslevagolemma}), 
we have a great freedom in specifying it what we now exploit as 
follows. Consider the original simply connected and closed $M$ used in Theorem 
\ref{fotetel1}. Simply connectedness implies the vanishing of the first de 
Rham cohomology of $M$ therefore if we put any Riemannian metric onto $M$ 
and consider the corresponding Laplacian on $1$-forms, its kernel is trivial. 
The Hodge decomposition theorem then says that any $1$-form $\xi$ on $M$ 
uniquely splits as $\xi=\dd f+\dd^*\eta$ where $f$ is a function and 
$\eta$ a $2$-form on $M$. The corresponding dual decomposition of a smooth 
vector field $v$ on $M$ therefore looks like $v={\rm grad}f+{\rm div}T$ 
where $T$ is a $(2,0)$-type tensor field. 

Motivated by this, consider now the space $X^\times$ of Theorem 
\ref{lorentztetel} and recall that it is homeomorphic to $M^\times$ 
consequently has vanishing first de Rham cohomology, too. Therefore, 
as a first and naive choice, we set the nowhere 
vanishing vector field (\ref{vektormezo}) used to construct the Ricci-flat 
Lorentzian metric $g_L$ on $X^\times$ out of the Ricci-flat Riemannian one 
$g$ to be of the form 
\begin{equation}
v:={\rm grad} f
\label{szetszedes} 
\end{equation}
where $f:X^\times\rightarrow (-\infty ,0]$ is a Morse function 
(to be defined shortly) on $X^\times$ such that $f^{-1}(-\infty)$ corresponds 
to the creased end of $X^\times$ while $f^{-1}(t)\subset 
X^\times$ are compact level sets for all $-\infty <t\leqq 0$ and in 
particular the point $f^{-1}(0)$ is the ``top'' of $X^\times$ 
(see Figure 6). 
\vspace{0.1in}

\centerline{
\begin{tikzpicture}[scale=0.6, rotate=-90]
\node at (15,-2) {$f:$};
\node at (15,3) {$\xrightarrow{\hspace{1cm}}$};
\node at (15,7) {$(-\infty ,0]$};
\draw [fill=black] (11.2,5) circle (0.15cm);
\draw [thick] (11.2,5) -- (18.15,5);
\draw [thick] (18.3,5) circle (0.15cm);
\node at (15,0) {$X^\times$};
\draw [thick] plot [tension=0.8,smooth] coordinates
{(11,0) (11.1,0.5) (11.5,1) (12.25,1.4) (13.5,1.5) (14.75,1.2) (16,0.2)};
\draw [thick] plot [tension=0.8,smooth] coordinates
{(11,0) (11.1,-0.5) (11.5,-1)
(12.25,-1.4) (13.5,-1.5) (14.75,-1.2) (16,-0.2)};
\draw [thick] plot [tension=0.8,smooth] coordinates {(16,0.2) (17,0.2) 
(18,1)};
\draw [thick] plot [tension=0.8,smooth] coordinates {(16,-0.2) (17,-0.2)
(18,-1)};
\draw [red, thick] plot coordinates {(18,1) (18.5,0.9) (18,0.7) (18.5,0.6)
(18,0.5) (18.5,0.4) (18,0.3) (18.5,0.2) (18,0.1) (18.5,0) (18,-0.1)
(18.5,-0.2) (18,-0.3) (18.5,-0.4) (18,-0.5) (18.5,-0.6) (18,-0.7)
(18.5,-0.9) (18,-1)};
\draw [thick] plot [tension=0.8,smooth] coordinates {(12.2,-0.8)
(11.9,0.1) (12.2,0.8)};
\draw [thick] plot [tension=0.8,smooth] coordinates {(12.1,-0.65)
(12.4,0.1)(12.1,0.7)};
%\label{abra6}
\end{tikzpicture}
}
\vspace{0.1in}

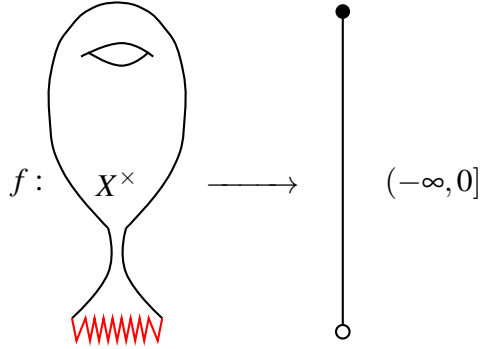
\centerline{Figure 6. The manifold $X^\times$ with a zig-zag 
representing its creased end} 
\centerline{and a Morse function $f:X^\times\rightarrow (-\infty,0]$ on it.}
\vspace{0.3in}

\noindent Moreover ${\rm grad} f$ in (\ref{szetszedes}) is defined by 
$\dd f =g({\rm grad}f\:,\:\cdot\:)$ to be the dual vector field of the 
$1$-form $\dd f$ with respect to the original {\it Riemannian} metric $g$ on 
$X^\times$. If the choice in (\ref{szetszedes}) is possible then we gain a 
very nice picture on the vacuum space-time $(X^\times ,g_L)$. Namely, 
${\rm grad}f :X^\times\rightarrow L\subset L\oplus L^\perp =TX^\times$ is a 
vector field such that for a generic $t\in (-\infty ,0]$ it does not 
vanish and the level set $f^{-1}(t)\subset X^\times$ is a $3$ dimensional 
closed (i.e., compact without boundary) submanifold with 
$Tf^{-1}(t)=L^\perp\subset L\oplus L^\perp=TX^\times$. Hence with 
respect to $g_L$ we find that ${\rm grad}f$ 
is a timelike and by definition future-directed vector field 
$g_L$-orthogonal for the level sets which are spacelike. In other words: 
{\it The vector field $v$ in (\ref{vektormezo}) is an infinitesimal 
observer in the space-time $(X^\times ,g_L)$. If $v$ has the form 
(\ref{szetszedes}) then $v$ can be identified with a global classical 
observer in the sense that the level value $t\in (-\infty ,0]$ 
corresponds to its global classical proper time as moves along its future 
directed own timelike curves (i.e., the integral curves of $v={\rm grad}f$)  
and the level sets $f^{-1}(t)\subset X^\times$ correspond to its global 
classical spacelike submanifolds}. However this picture is too 
naive because $f$ may attain critical points i.e., $p\in X^\times$ where 
${\rm grad}f(p)=0$ as we know from Morse theory. Hence the nowhere-vanishing 
vector field (\ref{vektormezo}) cannot globally look like (\ref{szetszedes}).

{\it A rapid course on Morse theory}. The following things are well 
known \cite{gom-sti,mil} but we summarize them here for completeness and 
convenience. Let $N$ be a smooth $n$-manifold. The point $p\in N$ is a 
{\it critical point} of a smooth function $f:N\rightarrow\R$ iff in a local 
coordinate system $(U_p, x_1,\dots,x_n)$ centered at $p$ all the 
partial derivatives vanish there i.e., $\partial_if(p)=0$ for all 
$i=1,\dots, n$ and it is {\it non-degenerate} iff the 
matrix $(\partial^2_{ij}f(p))_{i,j=1,\dots,n}$ is not singular. Moreover 
$c\in\R$ is a {\it critical value} iff the level set $f^{-1}(c)\subset N$ 
contains a critical point. The smooth function $f:N\rightarrow\R$ is a 
{\it Morse function} along $N$ iff it admits only non-degenerate 
critical points such that each critical value level set contains at most one 
critical point. (Being non-degenerate already implies that 
the critical points are isolated \cite[Corollary 2.3]{mil}). We shall also 
assume below that the 
level set $f^{-1}(c)\subset N$ is compact as well, for all $c\in\R$. 

We know the following things. If $c\in\R$ is non-critical 
then $f^{-1}(c)\subset N$ is a smooth $n-1$ dimensional submanifold. If 
$c\in\R$ critical with a single critical point $p\in f^{-1}(c)\subset N$ 
then (cf. \cite[Lemma 2.2]{mil}) there exists a local coordinate system 
$(U_p, y_1,\dots,y_n)$ about $p$ i.e., $y_1(p)=\dots=y_n(p)=0$, in which 
\[f\vert_{U_p}(y_1,\dots,y_n)=f(0,\dots,0)-\sum\limits_{i=1}^k
y^2_i+\sum\limits_{i=k+1}^{n}y^2_i\]
and the number $0\leqq k\leqq n$ is called the {\it index} of the critical 
point. Therefore a critical point of index $k=0$ is a local minimum while 
with index $k=n$ is a local maximum of $f$. Take $c\in\R$, $\varepsilon >0$ 
and suppose that $[c-\varepsilon,c+\varepsilon]\subset\R$ 
consists of non-critical values only. Then (cf. \cite[Theorem 3.1]{mil}) 
$f^{-1}(c-\varepsilon)$ and $f^{-1}(c+\varepsilon)$ are diffeomorphic. 
If the only critical value in $[c-\varepsilon,c+\varepsilon]$ 
is $c$ and its unique critical point $p\in f^{-1}(c)$ is of index $k$ then 
(cf. \cite[Theorem 3.2]{mil}) $f^{-1}(c+\varepsilon)$ is obtained from 
$f^{-1}(c-\varepsilon)$ by glueing to the boundary of 
$f^{-1}((-\infty ,c-\varepsilon ])$ a closed $n$-ball $B^n$ in the form of 
a $k$-handle $B^k\times B^{n-k}$. More precisely take an embedding 
$\varphi_k: S^{k-1}\times B^{n-k}\rightarrow f^{-1}(c-\varepsilon)$ 
and glue $B^n$ to $f^{-1}((-\infty, c-\varepsilon])$ by identifying 
\[S^{k-1}\times B^{n-k}\subseteqq\partial(B^k\times B^{n-k})=
(S^{k-1}\times B^{n-k})\cup (B^k\times S^{n-k-1})\] 
with the image $\varphi _k(S^{k-1}\times B^{n-k})\subseteqq
\partial\left(f^{-1}((-\infty, c-\varepsilon])\right)=f^{-1}(c-\varepsilon )$. Then 
after ``smoothing off the corners'' we obtain an $n$ dimensional 
manifold-with-boundary $f^{-1}((-\infty, c-\varepsilon])\cup_{\varphi_k}B^n$ 
and $f^{-1}(c+\varepsilon)$ is diffeomorphic to 
$\partial\left(f^{-1}((-\infty, c-\varepsilon])\cup_{\varphi_k}B^n\right)$. 
For instance if $k=0$ then $B^n$ is glued along $S^{-1}\times B^n$ where 
$S^{-1}=\emptyset$ i.e., it is not glued hence this critical point is a local 
minimum; while if $k=n$ then $B^n$ is attached along $S^{n-1}\times B^0$ 
where $B^0$ is a point i.e., it is attached along its full boundary 
$S^{n-1}$ hence this is a local maximum of $f$. Note that replacing the 
bottom-up function $f$ with the top-down function $-f$ critical points with 
index $k$ and $n-k$ interchange. 

Critical points necessarily occur. If $N$ is compact then a fundamental 
result of Morse theory (cf. \cite[Theorem 5.2]{mil}) states that if 
$m_k(N)$ denotes the number of critical points of index $k$ and $b_k(N)$ the 
$k^{{\rm th}}$ Betti number of $N$ then $b_k(N)\leqq m_k(N)<+\infty$. 
If $N$ is not compact then in general no such lower bounds exist but some 
$m_k(N)$'s can be even infinite. For further details cf. 
\cite[Chapter 4]{gom-sti} or \cite{mil}. 

Returning to our problem, we therefore correct (\ref{szetszedes}) as 
follows. Although critical points of $f$ are unavoidable, they are at least 
isolated i.e., for all $p,q\in X^\times$ pairs of critical points there exist 
small surrounding open neighbourhoods $U_p,U_q\subset X^\times$ such that 
$U_p\cap U_q=\emptyset$. Then taking the union  
\[C_f:=\bigcup\limits_{\mbox{$p$ is a critical point of $f$}}U_p\] 
which is therefore disjoint and supposing that this set is sharply 
concentrated around the critical 
points of $f$ in $X^\times$, let us correct (\ref{szetszedes}) to
\[v:={\rm grad}f+w\] 
where $w$ is a smooth vector field (of the form $w={\rm div}T$) on 
$X^\times$ such that $w(p)\not=0$ in the critical point $p$ but 
${\rm supp}\:w\subset C_f$ i.e., $w$ vanishes outside of 
$C_f\subset X^\times$. Fortunately this changes our physical picture 
on $(X^\times, g_L)$ only locally (i.e. close to a critical point 
only). More precisely, the classical observer picture of $v$ breaks down 
only in the vicinity of critical points of its Morse function part. Therefore 
from now on we suppose: {\it if $v={\rm grad}f+w$ is a non-vanishing 
vector field on $X^\times$ then the infinitesimal observer provided by $v$ 
in the original space-time $(X^\times, g_L)$ gives rise to a global classical 
observer in the sense above at least on the open domain
\begin{equation}
(X^\times\setminus\overline{C}_f\:,\: 
g_L\vert_{X^\times\setminus\overline{C}_f}\:)\subsetneqq (X^\times, g_L)
\label{megszoritas}
\end{equation}
because $v={\rm grad}f$ along this restriction as before}.

Let us ask ourselves now about the ``experiences'' of this partial global 
classical observer, constructed from a Morse function, as it moves in 
$(X^\times, g_L)$. That is, consider a Morse function $f$ on $X^\times$ as 
above (see Figure 6) with an associated global classical observer 
on the restricted domain $X^\times\setminus\overline{C}_f$. This observer has a 
global proper time $t\in (-\infty, 0]$ measured by $f$ with the 
infinite past $t=-\infty$ being the creased end of $X^\times$ and 
also has corresponding global spacelike 
$\Sigma_t:=f^{-1}(t)\subset X^\times\setminus\overline{C}_f$ for 
appropriate $t$'s which are closed $3$-manifolds. First, fix 
$-\infty < K<0$ such that $\Sigma_K$ is a submanifold and consider the 
compact part $f^{-1}([K, 0])\subsetneqq X^\times$. As the 
observer moves forwards in time i.e., from $t=K$ to $t=0$ along the integral 
curves of ${\rm grad}f$ then only finitely many critical points occur. As we 
have seen, around these points the spacelike $\Sigma_t$'s change topology by 
picking up a $k$-handle according to the index of the critical point. 

Now consider the much more interesting non-compact $f^{-1}((-\infty 
,K])\subset X^\times$ regime, the downward ``neck'' part in Figure 6. 
If $K<0$ is sufficiently small (we mean $\vert K\vert >0$ is 
sufficiently large) we can suppose that $f^{-1}((-\infty, K])$ is fully 
contained in the exotic but topologically trivial summand $R^4$ of $X^\times$ 
in its decomposition $X^\times =M\#R^4$. Therefore if $-\infty <t\leqq K$ 
then $\Sigma_t$ is fully contained in the $R^4$ summand. We can 
without loss of generality suppose that $\Sigma_K$ surrounds the 
attaching region of $M$ and $R^4$ hence $\Sigma_K$ is diffeomorphic to 
$S^3$. Now take an observer in $(X^\times,g_L)$ moving backwards in time along 
the integral curves of ${\rm grad} f$ i.e. from $t=K$ downwards $t=-\infty$. A 
generic value of $t$ 
is not critical for $f$ consequently the corresponding spacelike submanifold 
$\Sigma_t$ exists. Consider a fixed time $-\infty<t_0<K$ which is a 
critical value of $f$. How the corresponding transition between the 
$\Sigma_t$'s then looks like? As we have seen, in this moment always a single 
$4$-ball $B^4$, attached through its boundary $S^3$ in various ways to 
$\Sigma_t$ depending on the index $k$ of the critical point, is going to be 
removed from the latter space-time portion $f^{-1}([t_0,K])$. Therefore, as we 
move backwards in time provided by $f$ (or move forwards in time provided by 
$-f$) and pass through the moment $t_0$ the space $\Sigma_{t_0+\varepsilon}$ 
undergoes one of the following topological transitions: 
\begin{itemize}

%\item[$*$] If $k=0$ then at $t_0$ an $S^3$, disjoint from 
%$\Sigma_{t_0+\varepsilon}$, is annihilated (or equivalently, attached 
%completely, is created); 

\item[$*$] If $k=1$ then at $t_0$ an $S^3$, attached 
along two disjoint $B^3$'s to $\Sigma_{t_0+\varepsilon}$, is annihilated 
(or equivalently, attached along a thickened $S^2$, is created); 

\item[$*$] If $k=2$ then at $t_0$ an $S^3$, attached 
along a thickened knot to $\Sigma_{t_0+\varepsilon}$, is annihilated 
(or equivalently, attached along a thickened knot, is created); 

\item[$*$] If $k=3$ then at $t_0$ an $S^3$, attached 
along a thickened $S^2$ to $\Sigma_{t_0+\varepsilon}$, is annihilated 
(or equivalently, attached along two disjoint $B^3$'s, is created)

%\item[$*$] If $k=4$ then at $t_0$ an $S^3$, attached completely to 
%$\Sigma_{t_0+\varepsilon}$, is annihilated (or equivalently, a disjoint 
%one is created)

\end{itemize} 
and in this way the latter space $\Sigma_{t_0+\varepsilon}$ 
evolves into to the earlier $\Sigma_{t_0-\varepsilon}$ as we move backwards 
in time. Strictly mathematically speaking this $k$-handle attachment is to be 
performed ``instantaneosly'' somewhere along the singular level surface
$\Sigma_{t_0}$ carrying a unique critical point $p$ at the moment $t_0$; 
however from a physical viewpoint we can rather suppose that it occurs within 
the ``non-classical'' (with respect to the observer provided by ${\rm grad}f$) 
region $\Sigma_t\cap U_p\subset C_f$ at some unspecified time $t\in 
(t_0-\varepsilon ,t_0+\varepsilon)$ such that $\Sigma_{t_0\pm\varepsilon}
\cap U_p$ are still not empty (see Figure 7). Beside the $f$ Morse 
function picture, we have formulated all processes in the dual picture of 
the reversed Morse function $-f$ as well in order to gain full symmetry in 
the formulation. Moreover we note that applying diffeomorphisms on $X^\times$ 
(or equivalently, modifying $f$) we can assume that along $f^{-1}((-\infty 
,K])$ with $K<0$ the $k=0,4$ handle attachment steps corresponding to local 
minima and maxima do not occur.
\vspace{0.3in}

\centerline{
\begin{tikzpicture}[scale=0.7, rotate=-90]
\draw [thick, fill=gray] (7.5,2) circle (0.25cm);
\node at (8.2,2) {$\Sigma_{t_0-\varepsilon}\cap U_p$};
\draw [thick] (7.5,-1) circle (1cm);
\draw [thick] (7.5,5) circle (1cm);
\node at (7.5,9) {$\Sigma_{t_0-\varepsilon}$};
\draw [thick, fill=gray] (5,2) circle (0.7cm);
\draw [fill] (5,2) circle (0.1cm);
\node at (6.2,2) {$\Sigma_{t_0}\cap U_p$};
\node at (5,1.5) {$p$};
\def\samples{100}
\def\c{1.3}
\pgfmathsetmacro\cc{\c*\c}
\draw [thick, rotate=90, xshift=2cm, yshift=-5cm]
    (0,0) --
    plot[variable=\t,
         domain=-45+1/\samples:45-1/\samples,
         samples=\samples,
         smooth,]
    (\t:{\cc*sqrt(2*\cc*cos(2*(\t))})
    -- cycle;
\draw [thick, rotate=90, xshift=2cm, yshift=-5cm]
    (0,0) --
    plot[variable=\t,
         domain=180-45+1/\samples:180+45-1/\samples,
         samples=\samples,
         smooth,]
(\t:{\cc*sqrt(2*\cc*cos(2*\t))})
    -- cycle;
\node at (5,9) {$\Sigma_{t_0}$};
\draw [thick] (2,2) circle (1.6cm);
\draw [thick, fill=gray] (2,2) circle (0.25cm);
\node at (2.7,2) {$\Sigma_{t_0+\varepsilon}\cap U_p$};
\node at (2,9) {$\Sigma_{t_0+\varepsilon}$};
%\label{abra7}
\end{tikzpicture}
}
\vspace{0.3in}
         
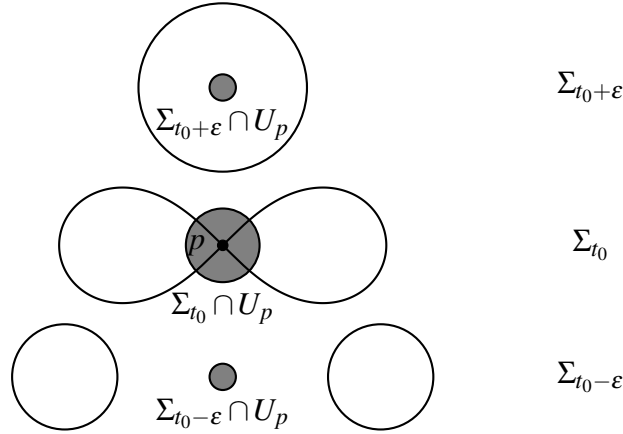
\centerline{Figure 7. Topology change about the critical point
$p\in\Sigma_{t_0}\cap U_p\subset X^\times$ about the moment 
$t_0\in (-\infty, 0]$.}
\vspace{0.3in}

Taking $-\infty\leftarrow t$ i.e., as moving backwards in time till the
creased end of $X^\times$ in Figure 6, in this process
the collection $\{\Sigma_t\}_{-\infty <t\leqq K}$ of spacelike
submanifolds looks like an evolution (in reversed
time) from $\Sigma_K=S^3$ into a three dimensional ``boiling foam'' limit
$\Sigma_{-\infty}$ or something like that. That is, these spacelike
submanifolds unboundedly continue to switch their topology; or in other
words the spatial oscillation between these states never stops and it is
reasonable to expect that all closed orientable $3$-manifolds
occur as $-\infty\leftarrow t$. Indeed, as
we noted in the Introduction, large exotic $R^4$'s always require
countably infinitely many handles in their handle decomposition therefore
moving backwards in time the $\Sigma_t$'s permanently continue changing their
topological type. Moreover soon or later $\Sigma_t$ very likely can be
arbitrary since the $k=2$ processes above are nothing but surgeries along
knots and all connected, closed, orientable $3$-manifolds arise this way
from $S^3=\Sigma_K$ by the Lickorish--Wallace theorem \cite{lic,wal}. This
``boiling foam'' picture therefore seems to be very weird and dynamical and
the sole ``driving force'' behind this dynamics is the non-standard
smooth structure along the end of $X^\times$. (Exactly the same thing is
responsible for the role of these spaces in $\overline{\bf SCCC}$, too.)
The existence of topologically different Cauchy surfaces in $\R^4$ is
already known to physicists, too \cite{new-cla}.

All the things have described up to this point might seem as mere 
mathematical nonsense. However they get even physically interesting if 
we recognize that this vivid spatial topology oscillation in $(X^\times, 
g_L)$ appears as a cosmological redshift phenomenon to our observer moving in 
(\ref{megszoritas}), as it looks back to the early creased end of $X^\times$ 
at late times. Let $E\in X^\times\setminus\overline{C}_f$ be a 
space-time event with a normalized future-directed timelike vector $n_E$ 
where a photon is emitted; in the geometrical optics approximation this 
photon travels along a future-directed null geodesic $\gamma$ in 
$(X^\times, g_L)$ till it is received in a later 
$R\in X^\times\setminus\overline{C}_f$ with corresponding receiver $n_R$. 
Taking any affine parameterization (i.e., $\nabla^L_{\gamma'}\gamma'=0$) the 
emitted frequency measured by $n_E$ is $\omega_E=-g_L(\gamma_E',n_E)$ while 
$\omega_R=-g_L(\gamma_R',n_R)$ is the frequency measured by the receiver. 
Then we define the {\it redshift factor} $z$ in the standard way by the 
frequency ratio 
\[1+z=\frac{\omega_E}{\omega_R}=\frac{g_L(\gamma_E',n_E)}
{g_L(\gamma_R',n_R)}\]
and say that the photon is {\it redshifted} along $\gamma$ if $z>0$. 
We adapt this general framework at least qualitatively to our setup as 
follows. Assume that the observer in the above process is given by 
$n=\frac{{\rm grad}f}{\vert{\rm grad}f\vert_{g_L}}$. Making use of the 
notation in the proof of Theorem \ref{lorentztetel}, 
${\rm grad} f$ is a section of $L\subset TX^\times$ hence 
$W_L{\rm grad}f=\sqrt{-1}\:{\rm grad}f$; moreover if 
$\gamma'=\gamma'_L+\gamma'_{L^\perp}$ is 
the unique decomposition according to $TX^\times =L\oplus L^\perp$ then 
$W_L\gamma' =W_L\gamma'_L+W_L\gamma'_{L^\perp}=
\sqrt{-1}\gamma'_L+\gamma'_{L^\perp}\in \sqrt{-1}\:L\oplus L^\perp$. 
Consequently 
\[g_L(\gamma'\:,\:n)=\frac{g_L(\gamma'\:,\:{\rm grad}f)}{\vert {\rm 
grad}f\vert_{g_L}}=\frac{g^\C(W_L\gamma'\:,\:W_L{\rm grad}f)}{\vert 
W_L{\rm grad}f\vert_{g^\C}}=\frac{-g(\gamma'_L\:,\:{\rm grad}f)}
{-\vert {\rm grad}f\vert_g}=\frac{g(\gamma'\:,\:{\rm grad}f)}{\vert 
{\rm grad}f\vert_g}\:\:.\]
Moreover 
\[\dd g(\gamma'\:,\:{\rm grad}f)\gamma'=-\dd g_L(\gamma'\:,\:{\rm 
grad}f)\gamma'=
-g_L(\nabla^L_{\gamma '}\gamma'\:,\:{\rm grad}f)-g_L(\gamma'\:,\:
\nabla^L_{\gamma '}{\rm grad}f)=-{\rm Hess}_f(\gamma',\gamma')\]
where ${\rm Hess}_f(x)=(\partial^2_{ij}f(x))_{i,j=1,\dots,4}$. Consider a 
non-critical point $q\in X^\times$ and its open neighbourhood $V_q\subset 
X^\times\setminus\overline{C}_f$ i.e. $V_q$ surely does not contain any 
critical point of $f$. Then there exists a local coordinate 
system $(V_q, t,x_1,x_2,x_3)$ centered at $q$ i.e. 
$t(q)=x_i(q)=0$ such that $f\vert_{V_q}(t,x_1,x_2,x_3)=t$ implying 
${\rm Hess}_f\vert_{V_q}=0$. Therefore $\dd g(\gamma'\:,\:{\rm grad}f)\gamma'=
\gamma'(g(\gamma'\:,\:{\rm grad}f))=0$ along $V_q$ i.e. if the photon path 
$\gamma$ does not intersect any critical point then $g(\gamma'\:,\:{\rm 
grad}f)$ is a non-zero constant along the whole $\gamma$. In this situation 
we end up with  
\[1+z=\frac{\vert {\rm grad}f(R)\vert_g}{\vert {\rm grad}f(E)\vert_g}\:\:.\]
As we emphasized throughout this note, the level surfaces 
$f^{-1}(t)\subset X^\times$ attain critical points more and more frequently 
as $-\infty\leftarrow t$. Consequently, the earlier space-time event 
$E\in f^{-1}(t_E)$ is ``more likely'' to be in the vicinity of a 
critical point $p_E\in f^{-1}(t_E)$ satisfying ${\rm grad}f(p_E)=0$ than 
the later event $R\in f^{-1}(t_R)$ with $t_R>t_E$. Therefore, acknowledging 
that a more careful statistical analysis is surely required, it is reasonable 
that ``typically'' $\vert{\rm grad}f(E)\vert_g\approx 0$ meanwhile 
$\vert{\rm grad} f(R)\vert_g\approx 1$ implying that the gradient ratio on 
the right hand side of $1+z$, when calculated for the ``typical'' early photon 
emitting event $E\in X^\times\setminus\overline{C}_f$ and late 
photon receiving event $R\in X^\times\setminus\overline{C}_f$, is large 
resulting in $z>0$. By the same reasoning this ratio even seems to be capable 
to be unbounded hence ``typically'' even $z>2$ seems reasonable which is 
exclusively characteristical for {\it cosmological} (i.e., not gravitational 
caused by a compact body, etc.) redshift. A cosmological context here is not 
surprising since our solutions $(X^\times ,g_L)$ are smooth while it has 
been known for a long time that in general relativity the gravitational field 
of an isolated massive object cannot be regular everywhere \cite{ein, ein-pau}.

Finally, one may raise the question about the place or role or relevance 
of this topology changing phenomenon within the full theory of 
(classical or even quantum) general relativity. Regarding 
this it is worth calling attention again that the Riemannian solutions 
$(X^\times ,g)$ underlying our smooth vacuum space-times $(X^\times 
,g_L)$ are not only Ricci-flat but even self-dual (see Theorem 
\ref{fotetel2} here). Consequently they are gravitational instantons and 
their appearance here looks reasonable for they are expected to generate 
these topology changes as tunnelings at the semi-classical (i.e. the leading 
term of quantum corrections) level. At first sight 
the whole picture presented here strongly resembles the structure of the 
vacuum sector of a non-Abelian gauge theory in temporal gauge over Minkowski 
space: in analogy with the present situation instantons of the Euclidean 
Yang--Mills theory over the Euclidean flat space execute semi-classical 
tunnelings between topologically (hence classically) separated classical 
vacua along space-like submanifolds in the original Minkowskian 
Yang--Mills theory over the Minkowskian flat space. 

However there is a subtle difference between the two tunneling processes 
which is probably worth recording here. In case of Yang--Mills theory 
all the aforementioned topologically different states connected by 
(anti)instanton effects are {\it vacua} hence the corresponding 
tunneling mechanism is time-symmetric which means that both instantons 
and antiinstantons occur and play a role. On the contrary in our 
gravitational situation the family $\{\Sigma_t\}_{-\infty <t\leqq 0}$ of 
topologically different spatial submanifolds with their corresponding 
Riemannian metrics inherited from their embeddings into $(X^\times,g_L)$ 
and connected by instanton effects are {\it not flat}; rather as 
$t\rightarrow 0$ this family looks like a sequence descending from 
quite complicated, topologically non-trivial highly curved compact $3$-spaces 
($\Sigma_t$'s with $t\ll 0$, the bottom part of Figure 6) 
towards topologically trivial $3$-spheres carrying metrics already close
to the standard round metric ($\Sigma_t$'s with $t\lessapprox 0$, 
the top of Figure 6). Therefore, as moving {\it forwards} in 
time the whole process seems to describe a sort of monotonic decay 
mechanism converting the gravitational degrees of freedom into other 
ones (like Yang--Mills fields, fermions, etc.) before reaching the 
gravitational vacuum (in our spatially compact situation the standard 
round $S^3$ plays the role of the flat geometry i.e. the gravitational 
vacuum). This process therefore seems to be not reversible and having a 
creased end introduces a sort of time direction along the cosmological 
space-time $(X^\times, g_L)$. Consequently the gravitational instantons 
provided by the spaces $(X^\times, g)$ are asymmetric unlike the 
gravitational instanton-antiinstanton pairs considered in \cite[Section 
III]{wit}. 

Are then $(X^\times, g)$'s physically relevant? Based on cluster 
decomposition Witten argues that a non-perturbative 
field is still relevant in a quantum theory if it is continuously deformable 
to the vacuum in an appropriate configuration space \cite[Section III]{wit}. 
Consider the case of traditional general 
relativity when space-time is topologically $\R^m$ and in particular the 
vacuum is the flat $\R^m$. Then by this argument gravitational 
instantons restricted to be exotic $m$-spheres if $m\not=4$. However if $m=4$ 
we cannot forget about exotic $\R^4$'s. In this exceptional situation we can 
follow Gompf \cite[Chapter 9.4]{gom-sti} 
and consider the configuration space $\crr_\sim$ of compact equivalence 
classes of smooth structures on $\R^4$. The set $\crr_\sim$ can be given the 
structure of a connected metrizable topological space with countable basis 
in which therefore the vacuum i.e. the standard $\R^4$ is represented by a 
point while our gravitational instanton $R^4$ by another point. 
Consequently within $\crr_\sim$ the gravitational instanton considered here 
is deformable into the vacuum. However the relevance of this purely formal 
observation is not clear neither from a physical nor a mathematical viewpoint.  
\vspace{0.1in}

\noindent {\bf Acknowledgement}. This paper is dedicated to Roger Penrose, the 
laureate of the 2020 Nobel Prize in physics. 
The author wishes to thank to I. K\'ad\'ar for the stimulating discussions and 
the Referee of JGP for raising lot of excellent clarifying questions.

\end{document}